\documentclass{article}

\usepackage{PRIMEarxiv}

\usepackage[utf8]{inputenc} % allow utf-8 input
\usepackage[T1]{fontenc}    % use 8-bit T1 fonts
\usepackage{hyperref}       % hyperlinks
\usepackage{url}            % simple URL typesetting
\usepackage{booktabs}       % professional-quality tables
\usepackage{amsfonts}       % blackboard math symbols
\usepackage{nicefrac}       % compact symbols for 1/2, etc.
\usepackage{microtype}      % microtypography
\usepackage{lipsum}
\usepackage{fancyhdr}       % header
\usepackage{graphicx}       % graphics

\usepackage{xcolor,colortbl}

\usepackage{boldline}% http://ctan.org/pkg/multirow
\usepackage{hhline}% http://ctan.org/pkg/hhline
\usepackage{amsmath}
\usepackage{amssymb}
\usepackage{algorithm}
\usepackage{algorithmic}
\usepackage{makecell}
\usepackage{adjustbox}
\usepackage{array, multirow, booktabs}

\DeclareMathOperator*{\argmax}{arg\,max}

\newcommand*\rot{\rotatebox{90}}

\makeatletter
\newcommand{\vast}{\bBigg@{4}}
\newcommand{\Vast}{\bBigg@{5}}
\makeatother

\graphicspath{{media/}}     % organize your images and other figures under media/ folder

\title{Risk-averse Stochastic Optimization for Farm Management Practices and Cultivar Selection Under Uncertainty
}

\author{
  Faezeh Akhavizadegan, Javad Ansarifar, Lizhi Wang \\
  Department of Industrial and Manufacturing Systems Engineering\\
  Iowa State University\\
  Ames, IA 50011, USA\\
  \texttt{\{faezeh, javad, lzwang\}@iastate.edu} \\
   \And
  Sotirios V. Archontoulis \\
  Department of Agronomy \\
  Iowa State University\\
  Ames, IA 50011, USA\\
  \texttt{sarchont@iastate.edu} \\
}

\begin{document}
\maketitle

\begin{abstract}
Optimizing management practices and selecting the best cultivar for planting play a significant role in increasing agricultural food production and decreasing environmental footprint. In this study, we develop optimization frameworks under uncertainty using conditional value-at-risk in the stochastic programming objective function. We integrate the crop model, APSIM, and a parallel Bayesian optimization algorithm to optimize the management practices and select the best cultivar at different levels of risk aversion. This approach integrates the power of optimization in determining the best decisions and crop model in simulating nature's output corresponding to various decisions. As a case study, we set up the crop model for 25 locations across the US Corn Belt. We optimized the management options (planting date, N fertilizer amount, fertilizing date, and plant density in the farm) and cultivar options (cultivars with different maturity days) three times: a) before, b) at planting and c) after a growing season with known weather. Results indicated that the proposed model produced meaningful connections between weather and optima decisions. Also, we found risk-tolerance farmers get more expected yield than risk-averse ones in wet and non-wet weathers.
\end{abstract}

% keywords can be removed
\keywords{Farm Optimization, Cultivar and Fertilizer Selection, Robust, Stochastic, Bayesian Optimization.}

\section{Introduction}

A key to improving agricultural output is developing precise decision-making in farming systems subject to different objectives and constraints \cite{dury2012models, ansarifar2020performance, ansarifar2021interaction}. This efficient farming system significantly impacts the farm's productivity and profitability in the short and long term. Although establishing this decision-making procedure plays a crucial role in assisting farmers, the development of such a decision-making process to consider nature behavior is a challenging issue \cite{ansarifar2020performance, ansarifar2022scheduling}. The challenge for optimizing the farming system is the dynamic and complex nature of the system as well as uncertainties in the future weather. Use of an optimization framework containing complex system processes and weather uncertainty to determine managerial decisions is a required step towards improving the management and efficiency in crop production systems \cite{ansarifar2021machine, akhavizadegan2021integration}. 

Among a variety of approaches, optimization-oriented approaches can provide a more promising and valuable decision-making platform for an agricultural system by taking account of agricultural, economic, and environmental factors at the same time \cite{garcia2017systems, russelle2007reconsidering}. Optimization techniques design optimal farming systems to achieve maximum returns and minimum environmental footprint corresponding to resource usage and management practices. In-depth reviews of operation research approaches in crop production systems have been published by Jain et al. \cite{jain2018optimization}, Caicedo Solano et al. \cite{caicedo2020towards}, Dury \cite{dury2012models}, Le Gal et al. \cite{le2011does}, Singh \cite{singh2012overview}, Martin et al. \cite{martin2013farming}, and Kanter et al. \cite{kanter2018evaluating}. The existing approaches were developed to address two types of problems: optimizing management strategies and optimizing land, water resources, and cropping patterns.

As the first category that focuses on optimization management strategies, Capitanescu et al. \cite{capitanescu2017multi} developed a multi-stage mixed-integer linear programming model for optimal farm management by considering environmental and crop rotation constraints. Klein et al. \cite{klein2013adapting}, by developing a regional multi-objective optimization approach, and Holzk{\"a}mper et al. \cite{ holzkamper2015assessing}, by proposing multi-objective optimization model, evaluated various management practices as a response to adaptation to climate change. To design a sustainable farming system, Liang et al. \cite{liang2018multi} formulated a multi-objective mixed-integer nonlinear fractional programming problem by considering the most important elements of organic farming systems including, crop cultivation, soil organic matter, nitrogen rate, and livestock production. Barbosa et al. \cite{barbosa2020modeling} proposed a convolutional neural network to assess yield response to nutrient and seed rate management. Groot et al. \cite{groot2012multi} designed a farming system using a multi-objective optimization model to maximize profit and organic matter accumulation and minimize labor hour and N loss as objective functions. Applications of dynamic programming models in optimizing agricultural management practices problems were conducted by Kennedy \cite{kennedy2012dynamic}, Parsons et al. \cite{parsons2009weed}, and Janov{\'a} et al. \cite{janova2014dynamic}. As the combinations of two categories, Pastori et al. \cite{pastori2017multi} developed a multi-objective evolutionary approach to assess the effects of water and nutrient management options in Africa. Moreover, Bhar et al. \cite{bhar2020coordinate} used the root zone water quality model and three global optimization methods to optimize the fertilization and irrigation decisions under the maximization of farm profit. Their study reported about 7\% improvement in yield and about 10\% improvement in profit compared to the experimental data of Trout and Bausch \cite{trout2017usda}.

The second category of existing research refers to the optimization of land, water resources, and cropping patterns. Optimization approaches for cropping pattern management include linear model \cite{singh2001optimal, alabdulkader2012optimization}, multi-stage linear programming model \cite{galan2015multi}, integer linear programming model \cite{dos2011crop, ansarifar2019new}, mixed-integer linear programming \cite{dogliotti2006influence}, multi-period mixed integer nonlinear programming model \cite{cervantes2020optimal}, stochastic programming \cite{vizvari2014stochastic}, fuzzy goal programming model \cite{ mirkarimi2013application, sharma2016fuzzy, soltani2011determining}, dynamic optimization model \cite{liu2016dynamic}, improved evolutionary model \cite{sarker2009improved}, pareto-based evolutionary model \cite{marquez2011multi}, crop rotation model \cite{schonhart2011croprota}, and network models \cite{salassi2013economically}. Most previous approaches that focus on optimizing land and water resources include linear programming model \cite{singh2012development, singh2012optimal, osama2017optimization, khare2007assessment}, multi-objective linear programming \cite{bartolini2007impact, galan2017multi}, mixed-integer nonlinear programming algorithms \cite{yue2013design}, multi-objective genetic algorithm \cite{lalehzari2016multiobjective, mousavi2017application}, multi-objective compromise programming \cite{cabrini2016modeling}, genetic algorithm and linear programming \cite{karamouz2009development, azamathulla2008comparison}, non-linear optimization model \cite{montazar2010conjunctive}, inexact rough-interval fuzzy linear programming \cite{lu2011inexact}, robust multistage interval-stochastic programming method \cite{li2011robust}, Simulation-optimization modeling \cite{singh2014simulation}, dynamic programming \cite{sarttra2013application}, chance-constrained linear programming \cite{singh2015land, sethi2006optimal}.

Despite significant previous efforts, our literature review on optimization of farm management reveals that the assessment of management practice options before and after gowning season with a realistic model that reflects the non-linear relationship between its explanatory variables is missing. Most previous studies used simple regression model to assess productivity impact. To our knowledge none has used a biophysical model. Hence, in this research, we used the APSIM crop model to mimic nature and simulate the farming system. To have a near-reality agricultural system, we used a calibrated crop model that simulates yields reasonably close to historical data. Then, this crop model is utilized as a simulator step of the optimization process to evaluate optimizers' decisions appropriately. 

In this paper, to address uncertainty in decision-making, we propose optimizers using parallel Bayesian algorithms as the core search engine to effectively search for the best combinations of management decisions over unknown objective functions. We incorporate risk considerations in the proposed decision-making framework via conditional value-at-risk (C-VaR) measure to make a decision under uncertainty and evaluate the the risk associated with the decisions during a growing season. This decision-making tool is able to take into account agricultural and farmers' objective functions (maximizing yield or profit, minimizing N loss and the environmental footprint), realistic farming system, and weather uncertainties to optimize the farm management practices before and after growing season and select the best management practice. In this research, farm management options are planting date, N fertilizer amount, fertilizing date, and plant density in the farm, and cultivar options are cultivars with different maturity days under uncertainty towards maximizing yield. We set up crop model from 25 locations across the US Corn Belt from 1984 to 2019 (34 years) and optimized management and cultivar options for three test years (2016-2018) at three time-wise strategies during growing season (first of March, first of May, and first of November). Finally, we interpret the optimal decisions of the proposed decision-making framework using stochastic and robust optimizers.

\section{Problem Definition}

The first crucial step of crop production for farmers is planning and decision making about which and when cultivars must be planted or what management practices must be carried out. The challenge for farmers during this decision-making process is weather uncertainties that determine yield potential and to some extend uncertainty in product prices. Farmers can not control the weather or the market price but can control the management (cultivar, fertilizer, and planting decisions). In this research, we consider a farmer who wants to make decisions about the management practices options and determines the best choice for the cultivar. The farmer seeks to optimize his objective function (maximizing the farm’s profit, minimizing N loss, and minimizing environmental footprint) for the next growing season by making the right decisions with unknown weather before the growing season. Management decision can be row spacing, plant density, planting date, fertilizing date, fertilizer amount, and type of cultivar. We aim to address this problem by making decision in three time in growing season such as first of March (earlier time), first of May (Plant time), and first of November (harvest time) and considering weather uncertainties simultaneously. The difference between the three time periods is the weather uncertainty.

\section{Method}

\begin{figure}[H]
	\centering
	\includegraphics[width=5in]{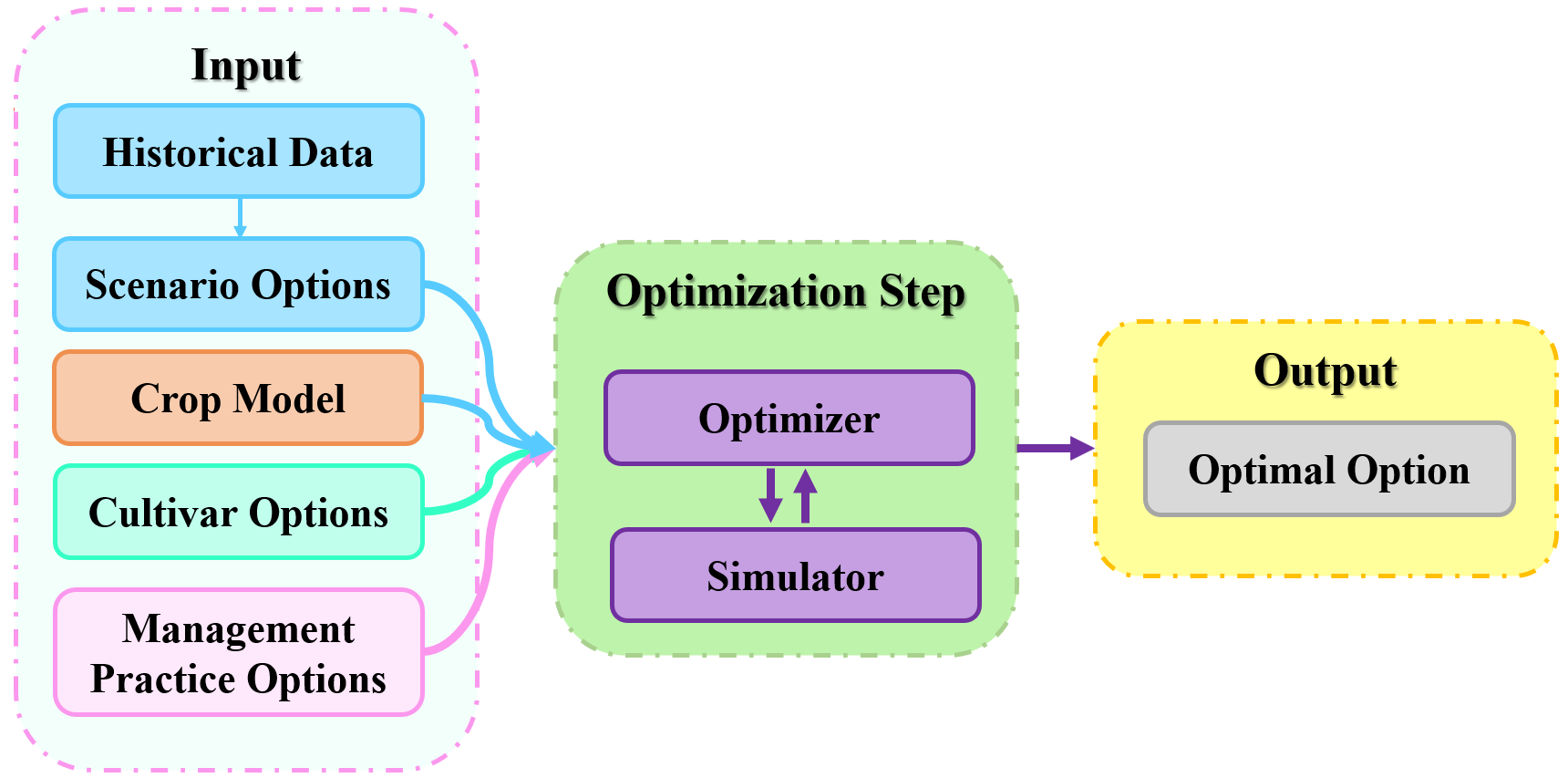}
	\caption{The proposed decision-making process for management practices and cultivars selection.}\label{prop}
\end{figure}

The proposed decision-making process for optimizing management practices and cultivars selection under uncertainty is illustrated in Figure \ref{prop}. Because of uncertain weather information and complex behavior of nature on management practices and cultivar options, we used the power of crop model and both stochastic and robust optimization frameworks. First, we defined different scenarios using historical information due to unknown weather. Then, the optimization model uses these scenarios to deal with a decision-making process under uncertainty. The optimization framework consists of two elements: optimizer and simulator. The framework is an iterative process because the optimizer determines the optimal solution during several iterations. Optimizer requires a simulator to evaluate the solution's performance. The optimizer is exploring the best solution in solution spaces, while the simulator's simulate an agriculture system based on each solution.

\noindent\textbf{Simulator}

We use the APSIM crop model \cite{holzworth2014apsim} as a near-reality agricultural system in our framework to evaluate the performance of specific management practices and selected cultivar. Crop models like Apsim can simulate nature processes by using a series of non-linear equations from the physiological understanding of plant processes \cite{heslot2014integrating, schauberger2017consistent}. Crop models have been successful in testing crop management strategies \cite{zhao2014sensitivity}, forecasting crop yield \cite{ basso2019seasonal}, evaluating cultivars in breeding programs \cite{lamsal2017efficient}, and exploring climate change impacts on productivity \cite{anwar2015climate}. The input of crop models is cultivar information, environment (soil profile and weather), and management practices. 

In the proposed model, the optimizer sends a decision to Apsim, and we modify the crop model's input based on the decision and run the model. According to Apsim output and evaluation criterion, the simulator returns the corresponding value.

\noindent\textbf{Optimizer}

As an optimizer, we propose a stochastic and robust optimization framework for optimizing farm management practices and cultivar options using a recently developed parallel Bayesian optimization algorithm \cite{akhavizadegan2021time}. We use several scenarios to cope with uncertainty (e.g., unknown weather data for the next growing season). In mathematical programming, we apply C-VaR as a tractable measure to model the risk associated with these scenarios. The core of the optimizer is a parallel Bayesian optimization algorithm. The proposed optimizer using a parallel Bayesian optimization is explained in Algorithm \ref{PBO}.

A parallel Bayesian optimization is an iterative approach so that it tries to improve the incumbent solution effectively. PBO runs N parallel BO \cite{mockus2012bayesian}, and at each iteration, BO models propose combinations of variables. Each BO proposes one combination of variables at each iteration and calls the crop model for all scenarios to evaluate its performance. In calculating an objective function, we use the C-VaR of yield performance to manage the risk associated with weather uncertainty. By defining the probability density function of yield under different scenarios and denoting $\alpha$ as a parameter indicating the right tail probability	of that function, C-VaR$_{\alpha}$ is defined as the expected value in the worst 100*$\alpha$\% of the yield distribution (yield set under different scenarios) \cite{rockafellar2000optimization}. For instance, if $\alpha=1$, the optimizer would be stochastic optimization. If $\alpha=1/|S|$ where $|S|$ is a number of scenarios and probability of scenarios are equal, the optimizer would be robust optimization. 

In a Stochastic approach, the objective function of one solution is the expected objective function values of all scenarios. In contrast, robust optimization computes it as the objective function value of the worst scenario. After calculating the objective function of one combination, we add this combination and associated objective function value to the dataset. Then, the gaussian process updates the posterior probability, and each BO proposes a new management combination. We continue until the maximum iteration is reached. In the end, The best combination of management practices and cultivar with the best objective function value between the dataset is chosen for the optimal solution.

\begin{algorithm}[H]
	\caption{\texttt{optimizer}} \label{PBO}
	\begin{algorithmic}[1]
		\STATE \textbf{Input:} Data set $D=\{ (x_0, y_0)\}$, $N$ as number BO models are run in parallel and $T$ as maximum iteration, $\alpha$ as a right tail probability of that yield distribution, $I_s$ as scenario $s \in S$, where $S$ denotes set of scenarios.
		\STATE \textbf{Output:} A local optimal solution $x^* \in \mathbb{R}^{1 \times p}$.
		\STATE \textbf{Step 0:} Set $K_n$, $\sigma$, and $u_n$ as a type of kernel, set of kernel parameters, and acquisition function of $n$th BO for all $n \in \{1,...,N\}$. Set incumbent solution $x^*=x_0$.
		\FOR{$t=1$ to $T$}
		\FOR{$n=1$ to $N$}
		\STATE \textbf{Step 1:} Use GP to update the posterior probability $\hat{f}_n$ and construct the acquisition function $u_n$.
		\STATE \textbf{Step 2:} Find $x_t^n$ by optimizing the acquisition function $u_n$ over function $f$: $x_t^n = \argmax_{x}u_n(x|D)$.
		\STATE \textbf{Step 3:} For each scenario $I_s$, call a crop model by modifying its file according to each scenario to evaluate the performance of solution $x_t^n$ for scenario $s$ as $f(x_t^n,I_s)$. 
		\STATE \textbf{Step 4:} Calculate the objective function of solution $x_t^n$ as mean of the $\alpha$-tail distribution of $f(x_t^n,I_s)$. For stochastic optimization as $y_t^n=\sum_{s \in S}{f(x_t^n,I_s)}/|S|$ and robust optimization as $y_t^n=min_{s \in S}\{f(x_t^n,I_s)\}$.
		\STATE \textbf{Step 5:} Augment the data $D = \{D,(x_t^n,y_t^n)\}$ and update the incumbent solution $x^*  = x_{\argmax (y)}$.
		\ENDFOR
		\ENDFOR
	\end{algorithmic}
\end{algorithm}

We provide additional details about the algorithm as follows.

\begin{itemize}
	
	\item \textbf{Hyperparameters:}

	Parameter $N$ denotes the number of parallel BO algorithms. By $T$, we define the maximum number of iterations as the stop criterion of the algorithm. This criterion can be replaced with a threshold so that the algorithm is stopped when the difference between the performance of incumbent and new solutions is less than a predefined threshold. $I_s$ denotes the information of scenario $s$. 
	
	\item \textbf{Step 0:} At this step, each BO chooses one type of acquisition function, kernel, and kernel parameters for its Gaussian process. 
	
	\item \textbf{Step 1:} In this step, Gaussian process updates the posterior probability, and build acquisition function.
	
	\item \textbf{Step 2:} The limited-memory quasi-Newton algorithm for bound-constrained optimization $\text{L-BFGS-B}$ \cite{zhu1997algorithm,byrd1995limited} in Python is applied to optimize the acquisition function. $\text{L-BFGS-B}$ algorithm is run with a set of random starting points to enhance the efficiency of optimizing the acquisition function. 
	
	\item \textbf{Step 3:} We modify the APSIM Files according to each new combination of variables from step 3 and run the APSIM for each scenario.
	
	\item \textbf{Step 4:} In this step, we calculate the objective function of the solution as mean of the $\alpha$-tail distribution of $f(x_t^n,I_s)$. If the optimization is Stochastic, the objective function value is the expected objective function values of all scenarios, where robust optimization calculates it as the objective function value of the worst scenario. 	
	
	\item \textbf{Step 5:} Then, the new solution and its objective function value are added to the dataset to expand BO's knowledge about posterior distribution. The incumbent solution is updated as the best combinations of parameters evaluated by APSIM under all scenarios in this step.
\end{itemize}

Indeed, stochastic optimization maximizes/minimizes the expected objective function value of all scenarios, while robust optimization tries to maximize/minimize the objective function value of the worst scenario. Both proposed stochastic and robust optimization frameworks ensure optimality and feasibility of solution even when the weather changes in the next growing season by balancing the level of uncertainty and objective function values. 

\section{Case Study}

To show the performance of the proposed framework in optimizing the managerial decision and selecting cultivar, we designed a case study 25 locations in the US Corn belt outlined by Akhavizadegan et al \cite{akhavizadegan2021time}. In each location (combinations of 5 counties $\times$ 5 soil series within a county) the model ran for 1984 to 2019 (34 years). The APSIM maize crop model \cite{keating2003overview} (version 7.9) was used because it is an open-source and widely used in the US Corn Belt and worldwide \cite{keating2018modelling, archontoulis2020predicting}. For each location, we used USDA-NASS management information to run the model (planting date, plant density, N-fertilizer) and used USDA-NASS corn yield data to calibrate cultivar parameters towards representing well actual yields. Overall, the model perform well in simulating corn yields (see Akhavizadegan et al \cite{akhavizadegan2021time}).

\subsection{Experiment Settings}

In this study, we considered 25 APSIM simulations corresponding to 25 locations, and we selected the following variables summarized in Table \ref{param} to optimize. The management options are planting date (10 possible values), N fertilizer amount (21 possible values), N fertilizing date (10 possible values), plant density (7 possible values), and cultivar options (4 possible values). All combination of variables resulted in 10*21*10*7 *4 = 58,800 simulations per weather scenario. Yield performance as objective functions was used for optimizing and evaluating the combination of management decisions and cultivar performance.

\begin{table}[H]
	\centering
	\caption{List of variables with their ranges and units used in this study.}\label{param}
	\begin{tabular}{ccc}
		\clineB{1-3}{2.5}
		\textbf{Variable} &\textbf{Value or Range}&\textbf{Unit} \\
		\hline
		Planting date& April 15th to June 15th (weekly steps)& day\\
		N fertilizer amount& 0 to 400 (20 kg/ha steps)& kg/ha \\
		N fertilizing date& Apr to June (weekly steps)& day\\
		Plant density& 2 to 14 (1 pl/m$^2$ plant steps)& pl/m$^2$\\
		Cultivars&100, 105, 110, and 115 relative maturity days&\\
		\clineB{1-3}{2.5}
	\end{tabular}
\end{table}

\begin{figure}[H]
	\centering
	\includegraphics[width=6in]{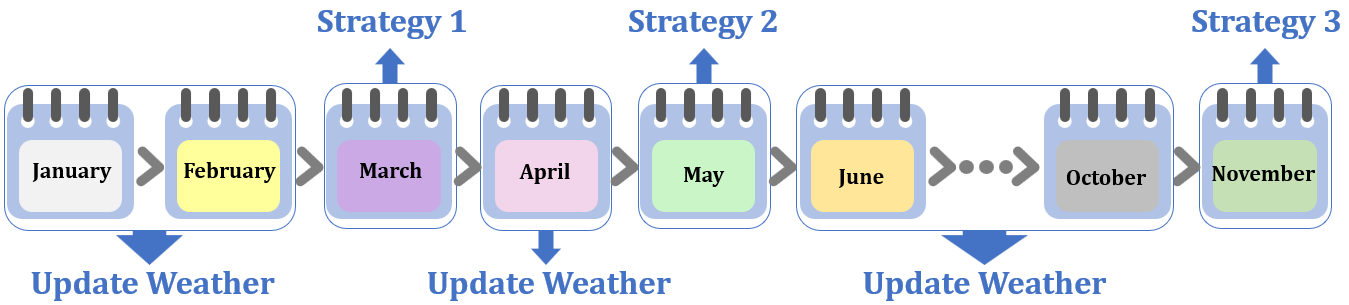}
	\caption{Overview of timewise strategies for decision.}\label{timely}
\end{figure}

Since farmers are choosing cultivars and deciding some management practices for the next growing season prior to planting, we explore the performance of the optimization framework under three timewise strategies during the growing season: before planting, at about planting, and at harvest. Figure \ref{timely} indicates the milestone of these strategies. In the first strategy, the decision is made on the first of March. The second decision is made on the first of May as the second strategy. Finally, we decide when all weather information is observed as the third strategy after crop harvesting (first of November), which will serve as a benchmark to compare decisions at the first and second decision times.  Hence, we have applied the optimization framework to decide on all 25 locations under all scenarios in three strategies. 

In this study, the weather is a source of uncertainties. Because the weather is part of the crop model's input, the unknown part of the weather data must be predicted before the growing season. Hence, weather scenarios are constructed using historical weather data to use as an unseen part of weather information. Since we have three strategies, we have to update each scenario's weather information by replacing it with observed weather information before decision time. The first strategy is the earliest time for making a decision, and the second one is close to the real, where the decision on the first of November is for comparing decisions under known weather with unknown weather.

We illustrate the effectiveness of our model by optimizing these variables for three test years (2016-2018) per location. The historical weather information from 1985 to year $Y-1$ is applied to generate the weather scenarios to optimize the variables at test year $Y$. Table \ref{weather} represents a summary of the precipitation amount of five locations in each selected county. Comparing the accumulated precipitation and its third percentile reveals that weather in 2018 was a wet year followed by 2016 and 2017. This insight helps us in investigating how weather affects optimal managerial decisions and cultivar selection. The proposed optimization framework was implemented in Python 3. We conducted this study on a workstation with a 3.4 GHz CPU and 16 GB memory. For each setting, we calculated the C-VaRa for $\alpha={1/|S|,0.25,0.5,0.75,1}$ where $|S|$ is the number of weather scenarios.

\begin{center}	
	\begin{table}[H]
		\centering \small
		\caption{Summary of the precipitation amount (in mm/year) of five locations in each selected county.}\label{weather}
		\begin{adjustbox}{width=1\textwidth}
			\begin{tabular}{c|ccc|ccc|ccc|ccc|ccc}
				\clineB{1-16}{4}
				\multirow{2}{*}{Metric}&\multicolumn{3}{c|}{Logan}&\multicolumn{3}{c|}{Greene}& \multicolumn{3}{c|}{Boone}& \multicolumn{3}{c|}{Keokuk}& \multicolumn{3}{c}{Obrien}\\
				\hhline{~---------------}
				&2016&2017&2018&2016&2017&2018&2016&2017&2018&2016&2017&2018&2016&2017&2018\\
				\clineB{1-16}{2}
				Mean&2.4&2.2& 3.1& 3.8& 2.9& 3.8& 2.6& 2.3& 3.3& 2.3& 2.0& 3.3& 2.5&2.0& 3.1\\
				Std&6.7& 6.6& 9.1& 9.4& 7.2& 10.5& 7.59& 6.3& 8.4& 6.6& 5.9& 9.7& 7.5&5.7& 9.6\\
				Sum&4483.4&4109.0&5819.1&7026.1&5359.5&6990.6&4800.5&4210.3&6050.5&4266.8&3698.7&6023.5&4724.5&3671.8&5833.0\\
				25\%&0.0&0.0&0.0& 0.0&0.0&0.0& 0.0&0.0&0.0& 0.0&0.0&0.0& 0.0&0.0&0.0\\
				50\%&0.0&0.0&0.0& 0.0&0.0&0.0& 0.0&0.0&0.0& 0.0&0.0&0.0& 0.0&0.0&0.0\\
				75\%&1.0& 0.6& 1.1& 2.1& 1.4& 2.2& 0.9& 0.7& 1.5& 0.9& 0.8& 1.5& 0.7& 0.8& 0.9\\
				Max&101.2& 73.2& 125.4& 89.3& 56.4& 124.2& 100.7& 52.9& 100.1& 74.7& 114.7& 125.8& 101.0&68.3& 103.6\\
				\clineB{1-16}{4}
			\end{tabular}
		\end{adjustbox}
	\end{table}	
\end{center}

\section{Results}

The results of the proposed framework on optimizing management practices and selecting cultivars for three test years (2016, 2017, and 2018) are summarized in Table \ref{com2016_2018}. The average estimated yields of five locations for five counties at the end of the growing season using optimal managerial decisions of different strategies reveal that observing more weather information at the initial stages of the growing season cannot guarantee better decisions. But decisions on the first of May result in better performance than decisions on the first of March when the weather is wet (the year 2018).

Our proposed model with the optimized practices was able to increase the estimated yields by 8.4\%, 9.7\%, and 5.4\% using optimal decisions at the first of March at three test years 2016, 2017, and 2018, respectively. Optimal decisions on May 1 improve the estimated yields by 9 \%, 14.7\%, and 8.6\% at three test years 2016, 2017, and 2018, respectively. For most counties, our framework can improve the original estimated yield of the crop model except Greene county in Indiana in 2016 and 2017. It is because of the unique weather profile of selected locations in this county in 2016 and 2017 that are different from their historical weather distribution. 

Performances of the proposed framework at the third strategy (model without any uncertainty) indicate improvement in yield from 2.1\% to 18.6\%. It also emphasizes that uncertain weather can play significant and unavoidable factors in farmers' decisions. Moreover, results reveal that risk-averse farmers (farmers with decisions based on the worst weather scenario) cannot make good decisions at the first of March and May at both wet year and non-wet years compared to risk-tolerance ones. While, risk-tolerance ones (farmers with decisions based on most weather scenarios) are successful in deciding optimal or close-to-optimal managerial practices and cultivars, especially in non-wet weather 2017.

\begin{center}	
	\begin{table}[H]
		\centering \small
		\caption{Observed yield in Bu/Acre, original crop model estimated yield in Bu/Acre, and average estimated yields (in Bu/Acre) of five locations for each selected county at the end of growing season 2016, 2017, and 2018 under optimal managerial decisions of different strategies. The optimizer tries to maximize the expected value in the worst 100*$\alpha$\% of the yield distribution, and in a Stochastic approach, the optimizer maximizes the expected yield of all scenarios. While, in a robust approach, the optimizer maximizes the yield of a worst-case scenario. The second column from the right contains the average estimated yields of five locations in a county under USDA and literature review decisions. The last column includes observed state-level yield. Percentages show the improvement of average estimated yields of five locations at the end of the growing season using optimal decisions compared with the original average estimated yield.}\label{com2016_2018}
		\begin{adjustbox}{width=1\textwidth}
			
			\begin{tabular}{ccc|ccccc|ccccc|c|c|c}
				\clineB{1-16}{4}
				\multirow{2}{*}{\rot{Year \hspace{1.5em} }}&\multirow{2}{*}{\rot{State \hspace{1.5em} }}&\multirow{2}{*}{\rot{County\hspace{1.25em}}}&\multicolumn{5}{c|}{Strategy 1}&\multicolumn{5}{c|}{Strategy 2}&\multirow{2}{*}{\rot{Strategy 3\hspace{0.9em}}}&\multirow{2}{*}{\rot{Est. Yield\hspace{0.9em}}}&\multirow{2}{*}{\rot{Obs. Yield\hspace{0.75em}}}\\
				\hhline{~~~----------~~~}
				&&&\rot{\hspace{0.25em}Robust}&\rot{$\alpha=0.25$}&\rot{\hspace{0.15em}$\alpha=0.5$}&\rot{$\alpha=0.75$}&\rot{Stochastic}&\rot{\hspace{0.25em}Robust}&\rot{$\alpha=0.25$}&\rot{\hspace{0.15em}$\alpha=0.5$}&\rot{$\alpha=0.75$}&\rot{Stochastic}&\multicolumn{1}{c|}{}&\multicolumn{1}{c|}{}\\
				\clineB{1-16}{2}	
				\multirow{10}{*}{\rot{2016}}&Illinois&Logan&\makecell{187.9\\-0.3\%}&\makecell{200.2\\6.2\%}&\makecell{202.4\\7.4\%}&\makecell{203.1\\7.8\%}&\cellcolor[gray]{0.85}\makecell{203.4\\7.9\%}&\makecell{192.9\\2.4\%}&\makecell{195.0\\3.4\%}&\cellcolor[gray]{0.85}\makecell{195.0\\3.4\%}&\makecell{190.8\\1.2\%}&\makecell{190.8\\1.2\%}&\makecell{214.0\\13.6\%}&188.5&219.4\\
				\clineB{2-16}{2}
				&Indiana&Greene&\makecell{169.5\\-9.8\%}&\makecell{178.5\\-5\%}&\cellcolor[gray]{0.85}\makecell{177.7\\-5.4\%}&\cellcolor[gray]{0.85}\makecell{177.7\\-5.4\%}&\makecell{177.0\\-5.7\%}&\makecell{174.0\\-7.4\%}&\cellcolor[gray]{0.85}\makecell{181.1\\-3.5\%}&\makecell{180.6\\-3.8\%}&\makecell{178.1\\-5.2\%}&\makecell{177.1\\-5.7\%}&\makecell{199.1\\6\%}&187.8&156\\
				\clineB{2-16}{2}
				&\multirow{6}{*}{Iowa}&Boone&\makecell{217.1\\4.9\%}&\makecell{222.2\\7.4\%}&\makecell{223.0\\7.8\%}&\makecell{223.0\\7.8\%}&\cellcolor[gray]{0.85}\makecell{224.3\\8.4\%}&\makecell{215.8\\4.3\%}&\makecell{220.8\\6.7\%}&\makecell{223.6\\8\%}&\makecell{223.3\\7.9\%}&\cellcolor[gray]{0.85}\makecell{225.5\\9\%}&\makecell{230.6\\11.5\%}&206.9&208.4\\
				\clineB{3-16}{2}
				&&Keokuk&\makecell{217.3\\4.6\%}&\makecell{219.1\\5.4\%}&\makecell{219.8\\5.8\%}&\cellcolor[gray]{0.85}\makecell{220.2\\6\%}&\makecell{218.7\\5.2\%}&\makecell{216.4\\4.2\%}&\makecell{218.8\\5.3\%}&\cellcolor[gray]{0.85}\makecell{218.9\\5.3\%}&\makecell{218.5\\5.1\%}&\makecell{217.9\\4.9\%}&\makecell{220.0\\5.9\%}&207.8&208.6\\
				\clineB{3-16}{2}
				&&Obrien&\makecell{257.5\\5.9\%}&\cellcolor[gray]{0.85}\makecell{259.8\\6.9\%}&\makecell{259.4\\6.7\%}&\makecell{259.4\\6.7\%}&\makecell{259.2\\6.6\%}&\makecell{259.7\\6.8\%}&\makecell{259.6\\6.8\%}&\makecell{259.2\\6.6\%}&\cellcolor[gray]{0.85}\makecell{259.8\\6.8\%}&\makecell{259.7\\6.8\%}&\makecell{265.1\\9\%}&243.1&210.7\\
				\clineB{1-16}{4}
				
				\multirow{10}{*}{\rot{2017}}&Illinois&Logan&\makecell{237.2\\-0.3\%}&\makecell{238.1\\0.1\%}&\makecell{245.2\\3\%}&\makecell{245.2\\3\%}&\cellcolor[gray]{0.85}\makecell{246.4\\3.5\%}&\makecell{234.9\\-1.3\%}&\makecell{238.8\\0.4\%}&\makecell{242.3\\1.8\%}&\makecell{244.2\\2.6\%}&\cellcolor[gray]{0.85}\makecell{246.4\\3.5\%}&\makecell{254.1\\6.8\%}&237.9&213\\
				\clineB{2-16}{2}
				&Indiana&Greene&\makecell{197.2\\-7.7\%}&\makecell{207.0\\-3.1\%}&\cellcolor[gray]{0.85}\makecell{210.8\\-1.3\%}&\makecell{209.6\\-1.9\%}&\makecell{209.0\\-2.2\%}&\makecell{202.3\\-5.3\%}&\makecell{204.9\\-4.1\%}&\makecell{205.3\\-3.9\%}&\makecell{205.3\\-3.9\%}&\cellcolor[gray]{0.85}\makecell{207.5\\-2.9\%}&\makecell{222.3\\4.1\%}&213.7&161.1\\
				\clineB{2-16}{2}
				&\multirow{6}{*}{Iowa}&Boone&\makecell{255.4\\2.2\%}&\cellcolor[gray]{0.85}\makecell{267.0\\6.8\%}&\makecell{266.4\\6.5\%}&\makecell{262.7\\5\%}&\makecell{265.3\\6.1\%}&\makecell{253.4\\1.3\%}&\makecell{262.3\\4.9\%}&\makecell{269.1\\7.6\%}&\makecell{269.4\\7.7\%}&\cellcolor[gray]{0.85}\makecell{271.2\\8.5\%}&\makecell{269.4\\7.7\%}&250.1&192.4\\
				\clineB{3-16}{2}
				&&Keokuk&\makecell{219.5\\-3.3\%}&\makecell{222.2\\-2.1\%}&\makecell{226.0\\-0.4\%}&\makecell{228.7\\0.8\%}&\cellcolor[gray]{0.85}\makecell{230.7\\1.6\%}&\makecell{219.9\\-3.1\%}&\makecell{222.6\\-1.9\%}&\makecell{223.9\\-1.4\%}&\makecell{226.2\\-0.3\%}&\cellcolor[gray]{0.85}\makecell{227.0\\0\%}&\makecell{231.7\\2.1\%}&227.0&176.4\\
				\clineB{3-16}{2}
				&&Obrien&\makecell{243.3\\1.8\%}&\makecell{246.4\\3\%}&\makecell{248.2\\3.8\%}&\makecell{255.7\\6.9\%}&\cellcolor[gray]{0.85}\makecell{262.4\\9.7\%}&\makecell{256.2\\7.1\%}&\makecell{265.8\\11.1\%}&\makecell{269.9\\12.9\%}&\makecell{270.1\\12.9\%}&\cellcolor[gray]{0.85}\makecell{274.2\\14.7\%}&\makecell{283.7\\18.6\%}&239.1&213.8\\
				\clineB{1-16}{4}
				
				\multirow{10}{*}{\rot{2018}}&Illinois&Logan&\makecell{214.6\\-1.1\%}&\makecell{215.1\\-0.9\%}&\cellcolor[gray]{0.85}\makecell{221.4\\2\%}&\cellcolor[gray]{0.85}\makecell{221.4\\2\%}&\makecell{221.3\\2\%}&\makecell{217.7\\0.3\%}&\makecell{225.4\\3.8\%}&\makecell{226.7\\4.4\%}&\cellcolor[gray]{0.85}\makecell{227.7\\4.9\%}&\makecell{227.5\\4.8\%}&\makecell{233.7\\7.7\%}&217.1&236.2\\
				\clineB{2-16}{2}
				&Indiana&Greene&\makecell{147.5\\-23.3\%}&\makecell{199.5\\3.8\%}&\makecell{198.0\\3\%}&\cellcolor[gray]{0.85}\makecell{203.1\\5.7\%}&\cellcolor[gray]{0.85}\makecell{203.1\\5.7\%}&\makecell{142.1\\-26.1\%}&\makecell{199.4\\3.8\%}&\makecell{202.1\\5.2\%}&\cellcolor[gray]{0.85}\makecell{208.6\\8.6\%}&\cellcolor[gray]{0.85}\makecell{208.6\\8.6\%}&\makecell{212.4\\10.5\%}&192.2&178.9\\
				\clineB{2-16}{2}
				&\multirow{6}{*}{Iowa}&Boone&\makecell{222.0\\2.9\%}&\makecell{223.7\\3.6\%}&\makecell{225.9\\4.7\%}&\cellcolor[gray]{0.85}\makecell{226.5\\4.9\%}&\cellcolor[gray]{0.85}\makecell{226.5\\4.9\%}&\makecell{223.2\\3.4\%}&\makecell{224.8\\4.2\%}&\makecell{225.8\\4.6\%}&\cellcolor[gray]{0.85}\makecell{226.1\\4.7\%}&\makecell{225.4\\4.4\%}&\makecell{239.4\\10.9\%}&215.8&193.7\\
				\clineB{3-16}{2}
				&&Keokuk&\makecell{215.4\\4\%}&\makecell{211.7\\2.2\%}&\cellcolor[gray]{0.85}\makecell{213.1\\2.9\%}&\cellcolor[gray]{0.85}\makecell{213.1\\2.9\%}&\makecell{210.8\\1.8\%}&\makecell{213.3\\3\%}&\makecell{216.2\\4.4\%}&\cellcolor[gray]{0.85}\makecell{216.8\\4.7\%}&\makecell{216.0\\4.3\%}&\makecell{214.6\\3.6\%}&\makecell{226.7\\9.5\%}&207.1&203.8\\
				\clineB{3-16}{2}
				&&Obrien&\makecell{239.8\\3.5\%}&\cellcolor[gray]{0.85}\makecell{240.2\\3.7\%}&\makecell{239.7\\3.5\%}&\makecell{238.5\\3\%}&\makecell{236.5\\2.1\%}&\makecell{240.7\\3.9\%}&\makecell{240.3\\3.7\%}&\cellcolor[gray]{0.85}\makecell{241.9\\4.4\%}&\cellcolor[gray]{0.85}\makecell{241.9\\4.4\%}&\makecell{240.7\\3.9\%}&\makecell{250.9\\8.3\%}&231.7&201.8\\
				\clineB{1-16}{4}
			\end{tabular}
		\end{adjustbox}
	\end{table}	
\end{center}

			\begin{center}	
				\begin{table}[H]
					\centering \small
					\caption{The percentage of time certain date was optimal N fertilizer date in three strategies. This percentage is respected to 25 different solutions (five levels of risk at five locations in each county). Three values in parentheses from left to right correspond to the first, second, and third strategies, respectively.}\label{NFD}
					\begin{adjustbox}{width=1\textwidth}
						\begin{tabular}{ccccccccccccc}
							\clineB{1-13}{4}
							\multirow{2}{*}{State}&\multirow{2}{*}{County}&\multirow{2}{*}{Year}&\multicolumn{10}{c}{N fertilizing date}\\
							\hhline{~~~----------}
							&&&1 Apr&8 Apr&15 Apr&22 Apr&29 Apr&6 May&13 May&20 May&27 May&3 Jun\\
							\clineB{1-13}{2}
							
							\multirow{3}{*}{Illinois}&\multirow{3}{*}{Logan}&2016&(0,40,40)&(0,60,0)&(20,0,20)&-&(20,0,0)&(20,0,20)&(40,0,20)&-&-&-\\
							&&2017&(40,20,0)&-&(20,40,0)&(0,20,0)&(20,20,20)&(20,0,40)&-&-&(0,0,20)&(0,0,20)\\
							&&2018&(40,60,0)&(0,20,0)&(0,0,40)&-&-&(40,0,0)&(0,0,20)&(20,20,40)&-&-\\
							\clineB{1-13}{2}
							\multirow{3}{*}{Indiana}&\multirow{3}{*}{Greene}&2016&-&-&-&-&-&-&(22.22,40,0)&(22.22,40,20)&(55.55,20,40)&(0,0,40)\\
							&&2017&-&-&-&-&-&-&(0,20,0)&(10,70,0)&(70,0,0)&(20,10,100)\\
							&&2018&-&-&-&-&-&-&-&-&(80,100,100)&(20,0,0)\\
							\clineB{1-13}{2}
							\multirow{9}{*}{Iowa}&\multirow{3}{*}{Boone}&2016& -&(20,0,0)&-&-&(0,20,0)&(20,20,20)&(20,20,20)&(20,0,20)&(0,20,20)&(20,20,20)\\
							&&2017&-&-&-&(10,20,0)&(20,40,0)&(30,0,0)&(10,0,0)&(30,0,0)&(0,20,20)&(0,20,80)\\
							&&2018&(10,20,0)&-&-&(10,0,0)&(0,0,20)&(0,20,20)&(20,20,0)&(20,20,40)&(20,0,20)&(20,20,0)\\
							\clineB{2-13}{2}
							&\multirow{3}{*}{Keokuk}&2016&-&-&-&-&-&-&-&-&(80,100,40)&(20,0,60)\\
							&&2017&-&-&(10,0,0)&(10,20,0)&-&(40,20,0)&(10,60,40)&-&(0,0,20)&(30,0,40)\\
							&&2018&-&(20,0,0)&-&(0,20,0)&-&(20,0,0)&-&(0,20,0)&(60,60,100)&-\\
							\clineB{2-13}{2}
							&\multirow{3}{*}{Obrien}&2016&-&(20,0,0)&-&(20,0,0)&(0,20,20)&(40,40,40)&(20,0,0)&(0,40,0)&(0,0,0,)&(0,0,40)\\
							&&2017&(20,10,0)&(0,10,0)&-&(20,0,0)&(20,20,0)&(40,30,0)&(0,10,0)&(0,20,20)&-&(0,0,80)\\
							&&2018&(0,20,20)&(20,0,0)&(20,0,20)&(20,0,20)&(0,0,20)&(20,50,0)&(20,0,20)&(0,20,0)&(0,10,0)&-\\
							\clineB{1-13}{4}
						\end{tabular}
					\end{adjustbox}
				\end{table}	
			\end{center}

Tables \ref{NFD}-\ref{pld} indicate the percentage of time certain decision variable's value was optimal in three strategies. This percentage is respected to 25 different solutions (five levels of risk at five locations in each county). According to Table \ref{NFD}, comparing fertilizing a farm in wet weather (2018) with non-wet weather (2017) reveals that it is better to fertilize in wet weather sooner than non-wet weather. Fertilizing at the end of May leads to a higher yield in Greene, Indiana, in wet and non-wet weathers. Optimal cultivars at different levels of risk aversion are shown in Table \ref{CUL}. It seems that cultivars with higher maturity dates are more favorable for most counties in three states in wet and non-wet weathers. Moreover, there is a relationship between the maturity date of cultivars and weather. In wet weather (2018), cultivars with lower maturity dates result in more yield. Table \ref{pde} shows that a higher plant density can guarantee a high yield at the end of the growing season, especially in wet weather. According to Table \ref{pld}, the proposed model suggests farmers plant in the middle of April for wet weather, where for non-wet weather, it finds the first two weeks of May as optimal planting dates.

					\begin{center}	
						\begin{table}[H]
							\centering \small
							\caption{The percentage of time certain cultivar was optimal cultivar in three strategies. This percentage is respected to 25 different solutions (five levels of risk at five locations in each county). Three values in parentheses from left to right correspond to the first, second, and third strategies, respectively.}\label{CUL}
							%		\begin{adjustbox}{width=1\textwidth}
								\begin{tabular}{ccccccc}
									\clineB{1-7}{4}
									\multirow{2}{*}{State}&\multirow{2}{*}{County}&\multirow{2}{*}{Year}&\multicolumn{4}{c}{Cultivar}\\
									\hhline{~~~----}
									&&& 100RM& 105RM& 110RM& 115RM\\
									\clineB{1-7}{2}
									\multirow{3}{*}{Illinois}&\multirow{3}{*}{Logan}&2016&(0,0,40)&(0,0,20)&-&(100,100,40)\\
									&&2017&-&-&(0,20,0)&(100,80,100)\\
									&&2018&-&-&-&(100,100,100)\\
									\clineB{1-7}{2}
									\multirow{3}{*}{Indiana}&\multirow{3}{*}{Greene}&2016&-&-&(0,0,40)&(100,100,60)\\
									&&2017&-&-&-&(100,100,100)\\
									&&2018&-&-&-&(100,100,100)\\
									\clineB{1-7}{2}
									\multirow{9}{*}{Iowa}&\multirow{3}{*}{Boone}&2016&-&-&-&(100,100,100)\\
									&&2017&-&-&-&(100,100,100)\\
									&&2018&-&(0,0,20)&(0,0,40)&(100,100,40)\\
									\clineB{1-7}{2}
									&\multirow{3}{*}{Keokuk}&2016&-&-&-&(100,100,100)\\
									&&2017&-&-&-&(100,100,100)\\
									&&2018&-&-&(0,0,100)&(100,100,0)\\
									\clineB{2-7}{2}
									&\multirow{3}{*}{Obrien}&2016&-&-&-&(100,100,100)\\
									&&2017&(20,0,0)&(20,0,0)&-&(60,100,100)\\
									&&2018&-&-&(0,0,60)&(100,100,40)\\
									\clineB{1-7}{4}
								\end{tabular}
								%\end{adjustbox}
							\end{table}	
						\end{center}

						\begin{center}	
							\begin{table}[H]
								\centering \small
								\caption{The percentage of time certain plant density was optimal plant density in three strategies. This percentage is respected to 25 different solutions (five levels of risk at five locations in each county). Three values in parentheses from left to right correspond to the first, second, and third strategies, respectively.}\label{pde}
								%		\begin{adjustbox}{width=1\textwidth}
									\begin{tabular}{cccccccccc}
										\clineB{1-10}{4}
										\multirow{2}{*}{State}&\multirow{2}{*}{County}&\multirow{2}{*}{Year}&\multicolumn{4}{c}{Plant density}\\
										\hhline{~~~-------}	
										&&&2&4&6&8&10&12&14\\
										\clineB{1-10}{2}
										\multirow{3}{*}{Illinois}&\multirow{3}{*}{Logan}&2016&-&-&-&(0,20,0)&(100,80,20)&(0,0,60)&(0,0,20)\\
										&&2017&-&-&-&(20,0,0)&(80,100,0)&(0,0,20)&(0,0,80)\\
										&&2018&-&-&-&(0,20,0)&(100,80,0)&-&(0,0,100)\\
										\clineB{1-10}{2}
										\multirow{3}{*}{Indiana}&\multirow{3}{*}{Greene}&2016&-&-&-&(44.44,60,60)&(55.55,40,40)&-&-\\
										&&2017&-&-&(0,10,0)&(30,0,0)&(70,90,20)&(0,0,40)&(0,0,40)\\
										&&2018&-&-&-&(20,40,0)&(80,60,0)&(0,0,40)&(0,0,60)\\
										\clineB{1-10}{2}
										\multirow{9}{*}{Iowa}&\multirow{3}{*}{Boone}&2016&-&-&-&-&(100,100,0)&(0,0,60)&(0,0,40)\\
										&&2017&-&-&-&(10,0,20)&(70,100,40)&(20,0,20)&(0,0,20)\\
										&&2018&-&-&-&-&(60,80,20)&(40,20,40)&(0,0,40)\\
										\clineB{2-10}{2}
										&\multirow{3}{*}{Keokuk}&2016&-&-&-&(20,0,0)&(80,100,40)&(0,0,40)&(0,0,20)\\
										&&2017&-&-&(0,0,40)&(40,20,60)&(50,60,0)&(10,20,0)&-\\
										&&2018&-&-&-&(60,20,20)&(20,60,0)&(20,20,60)&(0,0,20)\\
										\clineB{2-10}{2}
										&\multirow{3}{*}{Obrien}&2016&-&-&-&-&(100,100,0)&-&(0,0,100)\\
										&&2017&-&-&(20,0,0)&(20,0,0)&(40,90,0)&(0,10,0)&(20,0,100)\\
										&&2018&-&-&-&-&(80,100,0)&(20,0,0)&(0,0,100)\\
										\clineB{1-10}{4}
									\end{tabular}
									%\end{adjustbox}
								\end{table}	
							\end{center}

							\begin{center}	
								\begin{table}[H]
									\centering \small
									\caption{The percentage of time certain date was optimal planting date in three strategies. This percentage is respected to 25 different solutions (five levels of risk at five locations in each county). Three values in parentheses from left to right correspond to the first, second, and third strategies, respectively.}\label{pld}
									\begin{adjustbox}{width=1\textwidth}
										\begin{tabular}{ccccccccccccc}
											\clineB{1-13}{4}
											\multirow{2}{*}{State}&\multirow{2}{*}{County}&\multirow{2}{*}{Year}&\multicolumn{10}{c}{Planting date}\\
											\hhline{~~~----------}		
											
											&&&15 Apr &22 Apr &29 Apr &6 May &13 May &20 May &27 May &3 Jun &10 Jun &17 Jun\\
											\clineB{1-13}{2}
											\multirow{3}{*}{Illinois}&\multirow{3}{*}{Logan}&2016&(0,80,60)&(100,20,40)&-&-&-&-&-&-&-&-\\
											&&2017&(0,20,0)&(100,80,0)&-&(0,0,60)&(0,0,40)&-&-&-&-&-\\
											&&2018&(0,100,100)&(100,0,0)&-&-&-&-&-&-&-&-\\
											\clineB{1-13}{2}
											\multirow{3}{*}{Indiana}&\multirow{3}{*}{Greene}&2016&(100,100,40)&(0,0,60)&-&-&-&-&-&-&-&-\\
											&&2017&(100,90,0)&(0,10,0)&-&(0,0,60)&(0,0,40)&-&-&-&-&-\\
											&&2018&(80,80,80)&(20,20,20)&-&-&-&-&-&-&-&-\\
											\clineB{1-13}{2}
											\multirow{9}{*}{Iowa}&\multirow{3}{*}{Boone}&2016&(20,20,80)&-&(40,60,20)&(40,20,0)&-&-&-&-&-&-\\
											&&2017&(20,0,0)&(0,0,20)&(40,80,80)&(40,20,0)&-&-&-&-&-&-\\
											&&2018&(0,0,100)&-&(60,100,0)&(40,0,0)&-&-&-&-&-&-\\
											\clineB{2-13}{2}
											&\multirow{3}{*}{Keokuk}&2016&(0,0,20)&-&(20,0,0)&-&(60,60,0)&(20,40,60)&-&(0,0,20)&-&-\\
											&&2017&-&-&(30,0,20)&(0,0,20)&(40,20,20)&(30,80,40)&-&-&-&-\\
											&&2018&(0,0,100)&-&(40,0,0)&(20,0,0)&(20,60,0)&(20,40,0)&-&-&-&-\\
											\clineB{2-13}{2}
											&\multirow{3}{*}{Obrien}&2016&(0,20,20)&(80,60,0)&(20,20,80)&-&-&-&-&-&-&-\\
											&&2017&(20,20,100)&(40,60,0)&(0,20,0)&-&-&(20,0,0)&-&-&(20,0,0)&-\\
											&&2018&(0,20,80)&(80,60,0)&(20,20,20)&-&-&-&-&-&-&-\\

											\clineB{1-13}{4}
										\end{tabular}
									\end{adjustbox}
								\end{table}	
							\end{center}

\section{Conclusions}
This paper introduces a risk-averse stochastic model to optimize farm management practices and select cultivars on different levels of risk aversion and weather uncertainty. The proposed optimization framework consists of a crop model as a simulator and a parallel Bayesian optimization algorithm as an optimizer. Our framework can address farm management optimization issues under uncertainty at the same time. A computational experiment using 25 environments was conducted to illustrate the performance of the proposed framework. Five variables, including planting date, N fertilizer amount, N fertilizing date, plant density, and cultivars, were optimized to maximize yield at three time periods (first of March, May, and November). The optimization results produced insights about the impact of weather on yield by revealing that observing more weather information (between March to May) cannot guarantee better decisions in May than in March. The proposed model can boost the yield performance at the first of March and May compared with the original estimated yields of the crop model. Results indicate that risk-tolerance farmers can harvest more expected yield than risk-averse farmers. Moreover, by exploring the connection between weather and optimal decisions, the results reported meaningful agronomic insights between them.

\section{Conflict of Interest Statement}
The authors declare that the research was conducted in the absence of any commercial or financial relationships that could be construed as a potential conflict of interest.

\section*{Funding}
This work was partially supported by the National Science Foundation under the LEAP HI and GOALI programs (grant number 1830478) and under the EAGER program (grant number 1842097) and by the Plant Sciences Institute at Iowa State University.

\section{Author Contributions}
F.A., J.A., L.W, and S.V. designed the research questions. F.A., and J.A. prepared and cleaned the database, performed the experiment, statistical analysis, and analyzed the data set. F.A. designed and implemented a new algorithm, and created the figures. F.A., J.A, L.W., and S.V. interpreted  experiment results. F.A. wrote the manuscript. L.W. and S.V. thoroughly reviewed the manuscript and provided feedback. All authors read and approved the final manuscript.

%Bibliography
\bibliographystyle{unsrt}  
\bibliography{references}

	\newpage
\section*{Appendix}

\begin{center}	
	\begin{table}[H]
		\centering \small
		\caption{Observed yield, original crop model estimated yield, and estimated yields of five locations at each county at the end of growing season 2016 under optimal managerial decisions of different levels of risk aversion. The last column includes observed state-level yield. Percentages show the improvement of estimated yields of five locations at the end of the growing season using optimal decisions than the original estimated yield. }\label{2016ap}
		\begin{adjustbox}{width=1\textwidth}
			
			\begin{tabular}{cc|ccccc|ccccc|c|c|c}
				\clineB{1-15}{4}
				\multirow{2}{*}{\rot{State \hspace{1.5em} }}&\multirow{2}{*}{\rot{Location\hspace{1.25em}}}&\multicolumn{5}{c|}{Strategy 1}&\multicolumn{5}{c|}{Strategy 2}&\multirow{2}{*}{\rot{Strategy 3\hspace{0.9em}}}&\multirow{2}{*}{\rot{Est. Yield\hspace{0.9em}}}&\multirow{2}{*}{\rot{Obs. Yield\hspace{0.75em}}}\\
				\hhline{~~----------~~~}
				&&\rot{\hspace{0.25em}Robust}&\rot{$\alpha=0.25$}&\rot{\hspace{0.15em}$\alpha=0.5$}&\rot{$\alpha=0.75$}&\rot{Stochastic}&\rot{\hspace{0.25em}Robust}&\rot{$\alpha=0.25$}&\rot{\hspace{0.15em}$\alpha=0.5$}&\rot{$\alpha=0.75$}&\rot{Stochastic}&\multicolumn{1}{c|}{}&\multicolumn{1}{c|}{}\\
				\clineB{1-15}{3}			
				\multirow{10}{*}{Illinois}&logan-1&\makecell{213.8\\5.6\%}&\makecell{221\\9.2\%}&\makecell{221\\9.2\%}&\makecell{224.5\\10.9\%}&\makecell{224.5\\10.9\%}&\makecell{221\\9.2\%}&\makecell{221\\9.2\%}&\makecell{221\\9.2\%}&\makecell{200.2\\-1.1\%}&\makecell{200.2\\-1.1\%}&\makecell{229.2\\13.2\%}&202.4&\multirow{10}{*}{219.4}\\
				\hhline{~-------------~}
				&logan-2&\makecell{177.7\\1.1\%}&\makecell{177.7\\1.1\%}&\makecell{185.4\\5.5\%}&\makecell{182.3\\3.7\%}&\makecell{182.3\\3.7\%}&\makecell{167.3\\-4.9\%}&\makecell{182.5\\3.8\%}&\makecell{182.5\\3.8\%}&\makecell{182.5\\3.8\%}&\makecell{182.5\\3.8\%}&\makecell{206.2\\17.3\%}&175.79&\\			\hhline{~-------------~}
				&logan-3&\makecell{173.5\\-4.6\%}&\makecell{190.8\\4.8\%}&\makecell{190.8\\4.8\%}&\makecell{190.8\\4.8\%}&\makecell{192.3\\5.7\%}&\makecell{183.6\\0.9\%}&\makecell{190.3\\4.6\%}&\makecell{190.3\\4.6\%}&\makecell{190.3\\4.6\%}&\makecell{190.3\\4.6\%}&\makecell{207.8\\14.2\%}&181.98&\\
				\hhline{~-------------~}
				&logan-4&\makecell{180.5\\-2.5\%}&\makecell{199.4\\7.7\%}&\makecell{202.6\\9.5\%}&\makecell{202.6\\9.5\%}&\makecell{202.6\\9.5\%}&\makecell{180.5\\-2.5\%}&\makecell{190\\2.7\%}&\makecell{190\\2.7\%}&\makecell{190\\2.7\%}&\makecell{190\\2.7\%}&\makecell{209.5\\13.2\%}&185.11&\\
				\hhline{~-------------~}
				&logan-5&\makecell{194\\-1.7\%}&\makecell{212.4\\7.6\%}&\makecell{212.4\\7.6\%}&\makecell{215.5\\9.2\%}&\makecell{215.5\\9.2\%}&\makecell{212.4\\7.7\%}&\makecell{191.2\\-3.1\%}&\makecell{191.2\\-3.1\%}&\makecell{191.2\\-3.1\%}&\makecell{191.2\\-3.1\%}&\makecell{217.8\\10.4\%}&197.3&\\
				\clineB{1-15}{3}
				\multirow{10}{*}{Indiana}&greene-1&\makecell{177.3\\-6.5\%}&\makecell{176.5\\-6.9\%}&\makecell{170.2\\-10.2\%}&\makecell{170.2\\-10.2\%}&\makecell{170.2\\-10.2\%}&\makecell{168.4\\-11.1\%}&\makecell{179.2\\-5.5\%}&\makecell{178.8\\-5.7\%}&\makecell{170.2\\-10.2\%}&\makecell{170.2\\-10.2\%}&\makecell{201\\6.1\%}&189.51&\multirow{10}{*}{156}\\
				\hhline{~-------------~}
				&greene-2&\makecell{184\\-4.1\%}&\makecell{184\\-4.1\%}&\makecell{186.3\\-3\%}&\makecell{186.3\\-3\%}&\makecell{182.9\\-4.8\%}&\makecell{184\\-4.1\%}&\makecell{188\\-2.1\%}&\makecell{188\\-2.1\%}&\makecell{188\\-2.1\%}&\makecell{182.9\\-4.8\%}&\makecell{199.9\\4.1\%}&191.98&\\
				\hhline{~-------------~}
				&greene-3&\makecell{168.9\\-9.3\%}&\makecell{176.3\\-5.3\%}&\makecell{176.3\\-5.3\%}&\makecell{176.3\\-5.3\%}&\makecell{176.3\\-5.3\%}&\makecell{168.9\\-9.3\%}&\makecell{177.6\\-4.6\%}&\makecell{177.6\\-4.6\%}&\makecell{177.6\\-4.6\%}&\makecell{177.6\\-4.6\%}&\makecell{196.5\\5.6\%}&186.12&\\
				\hhline{~-------------~}
				&greene-4&\makecell{165.7\\-12.2\%}&\makecell{172.3\\-8.7\%}&\makecell{172.3\\-8.7\%}&\makecell{172.3\\-8.7\%}&\makecell{172.3\\-8.7\%}&\makecell{184\\-2.5\%}&\makecell{178.6\\-5.3\%}&\makecell{176.5\\-6.4\%}&\makecell{172.3\\-8.7\%}&\makecell{172.3\\-8.7\%}&\makecell{202.4\\7.3\%}&188.66&\\
				\hhline{~-------------~}
				&greene-5&\makecell{151.8\\-17.1\%}&\makecell{183.7\\0.3\%}&\makecell{183.7\\0.3\%}&\makecell{183.7\\0.3\%}&\makecell{183.7\\0.3\%}&\makecell{165\\-9.9\%}&\makecell{182.6\\-0.3\%}&\makecell{182.6\\-0.3\%}&\makecell{182.6\\-0.3\%}&\makecell{182.6\\-0.3\%}&\makecell{195.9\\7\%}&183.05&\\
				\clineB{1-15}{3}
				\multirow{30}{*}{Iowa}&boone-1&\makecell{202.8\\10.5\%}&\makecell{196.1\\6.9\%}&\makecell{198.7\\8.3\%}&\makecell{198.7\\8.3\%}&\makecell{205.3\\11.9\%}&\makecell{202.9\\10.6\%}&\makecell{196.1\\6.9\%}&\makecell{199.5\\8.7\%}&\makecell{198.2\\8\%}&\makecell{209.3\\14.1\%}&\makecell{215.3\\17.3\%}&183.49&\multirow{10}{*}{208.4}\\
				\hhline{~-------------~}
				&boone-2&\makecell{221.7\\2.9\%}&\makecell{231.5\\7.4\%}&\makecell{231.5\\7.4\%}&\makecell{231.5\\7.4\%}&\makecell{231.5\\7.4\%}&\makecell{230.2\\6.8\%}&\makecell{228.7\\6.1\%}&\makecell{231.5\\7.4\%}&\makecell{231.5\\7.4\%}&\makecell{231.5\\7.4\%}&\makecell{239.1\\11\%}&215.47&\\
				\hhline{~-------------~}
				&boone-3&\makecell{230.1\\5.3\%}&\makecell{230.1\\5.3\%}&\makecell{230.1\\5.3\%}&\makecell{230.1\\5.3\%}&\makecell{230.1\\5.3\%}&\makecell{227.3\\4.1\%}&\makecell{227.3\\4.1\%}&\makecell{232.1\\6.2\%}&\makecell{232.1\\6.2\%}&\makecell{232.1\\6.2\%}&\makecell{235.4\\7.8\%}&218.43&\\
				\hhline{~-------------~}
				&boone-4&\makecell{214.5\\3.3\%}&\makecell{226.6\\9.1\%}&\makecell{228.2\\9.9\%}&\makecell{228.2\\9.9\%}&\makecell{228.2\\9.9\%}&\makecell{218.1\\5\%}&\makecell{225.4\\8.5\%}&\makecell{228.2\\9.9\%}&\makecell{228.2\\9.9\%}&\makecell{228.2\\9.9\%}&\makecell{230.1\\10.8\%}&207.68&\\
				\hhline{~-------------~}
				&boone-5&\makecell{216.8\\3.3\%}&\makecell{226.8\\8.1\%}&\makecell{226.8\\8.1\%}&\makecell{226.8\\8.1\%}&\makecell{226.8\\8.1\%}&\makecell{200.7\\-4.4\%}&\makecell{226.8\\8.1\%}&\makecell{226.8\\8.1\%}&\makecell{226.8\\8.1\%}&\makecell{226.8\\8.1\%}&\makecell{233.5\\11.3\%}&209.86&\\
				\clineB{2-15}{3}
				&keokuk-1&\makecell{219.2\\2.6\%}&\makecell{220.5\\3.3\%}&\makecell{228.5\\7\%}&\makecell{228.5\\7\%}&\makecell{228.5\\7\%}&\makecell{218.2\\2.2\%}&\makecell{220.5\\3.2\%}&\makecell{220.5\\3.2\%}&\makecell{220.5\\3.2\%}&\makecell{220.5\\3.2\%}&\makecell{234.3\\9.7\%}&213.52&\multirow{10}{*}{208.6}\\
				\hhline{~-------------~}
				&keokuk-2&\makecell{219.2\\5.1\%}&\makecell{219.2\\5.1\%}&\makecell{216.7\\4\%}&\makecell{216.7\\4\%}&\makecell{218.8\\5\%}&\makecell{218\\4.6\%}&\makecell{218.1\\4.6\%}&\makecell{218.1\\4.6\%}&\makecell{218.1\\4.6\%}&\makecell{215.3\\3.3\%}&\makecell{212.7\\2\%}&208.5&\\
				\hhline{~-------------~}
				&keokuk-3&\makecell{215.6\\3.3\%}&\makecell{221.2\\5.9\%}&\makecell{221.2\\5.9\%}&\makecell{221.2\\5.9\%}&\makecell{219.6\\5.2\%}&\makecell{215.6\\3.3\%}&\makecell{221.1\\5.9\%}&\makecell{221.2\\5.9\%}&\makecell{218.6\\4.7\%}&\makecell{218.6\\4.7\%}&\makecell{218.1\\4.5\%}&208.78&\\
				\hhline{~-------------~}
				&keokuk-4&\makecell{215.5\\8.1\%}&\makecell{217.6\\9.1\%}&\makecell{215.5\\8.1\%}&\makecell{217.6\\9.1\%}&\makecell{209.3\\5\%}&\makecell{215.9\\8.3\%}&\makecell{217.2\\9\%}&\makecell{217.6\\9.1\%}&\makecell{218.1\\9.4\%}&\makecell{218.1\\9.4\%}&\makecell{217.9\\9.3\%}&199.36&\\
				\hhline{~-------------~}
				&keokuk-5&\makecell{217.4\\4\%}&\makecell{217.4\\4\%}&\makecell{217.4\\4\%}&\makecell{217.4\\4\%}&\makecell{217.4\\4\%}&\makecell{214.7\\2.7\%}&\makecell{217.4\\4\%}&\makecell{217.4\\4\%}&\makecell{217.4\\4\%}&\makecell{217.4\\4\%}&\makecell{217.4\\4\%}&208.99&\\
				\clineB{2-15}{3}
				&obrien-1&\makecell{259.8\\4.7\%}&\makecell{261.2\\5.3\%}&\makecell{261.2\\5.3\%}&\makecell{261.2\\5.3\%}&\makecell{260.1\\4.8\%}&\makecell{262.4\\5.7\%}&\makecell{260.3\\4.9\%}&\makecell{260.3\\4.9\%}&\makecell{260.8\\5.1\%}&\makecell{260.8\\5.1\%}&\makecell{265.9\\7.1\%}&248.18&\multirow{10}{*}{210.7}\\
				\hhline{~-------------~}
				&obrien-2&\makecell{260.2\\6.7\%}&\makecell{264.1\\8.3\%}&\makecell{262.1\\7.5\%}&\makecell{262.1\\7.5\%}&\makecell{262.1\\7.5\%}&\makecell{261.6\\7.3\%}&\makecell{262.1\\7.5\%}&\makecell{262.1\\7.5\%}&\makecell{264.1\\8.3\%}&\makecell{262.8\\7.8\%}&\makecell{268.1\\9.9\%}&243.91&\\
				\hhline{~-------------~}
				&obrien-3&\makecell{258.7\\6.8\%}&\makecell{259.7\\7.2\%}&\makecell{259.7\\7.2\%}&\makecell{260.8\\7.7\%}&\makecell{260.8\\7.7\%}&\makecell{258.7\\6.8\%}&\makecell{259.7\\7.2\%}&\makecell{259.7\\7.2\%}&\makecell{259.7\\7.2\%}&\makecell{260.7\\7.6\%}&\makecell{265.2\\9.5\%}&242.28&\\
				\hhline{~-------------~}
				&obrien-4&\makecell{254.5\\4.4\%}&\makecell{259.5\\6.5\%}&\makecell{259.5\\6.5\%}&\makecell{258.2\\6\%}&\makecell{258.2\\6\%}&\makecell{259.3\\6.4\%}&\makecell{259.5\\6.5\%}&\makecell{259.5\\6.5\%}&\makecell{259.5\\6.5\%}&\makecell{259.5\\6.5\%}&\makecell{264.8\\8.6\%}&243.73&\\
				\hhline{~-------------~}
				&obrien-5&\makecell{254.8\\7.2\%}&\makecell{254.8\\7.2\%}&\makecell{254.8\\7.2\%}&\makecell{254.8\\7.2\%}&\makecell{254.8\\7.2\%}&\makecell{256.7\\8\%}&\makecell{256.7\\8\%}&\makecell{254.6\\7.1\%}&\makecell{255\\7.2\%}&\makecell{255\\7.2\%}&\makecell{261.9\\10.1\%}&237.77&\\
				\clineB{1-15}{4}
			\end{tabular}
		\end{adjustbox}
	\end{table}	
\end{center}

\begin{center}	
	\begin{table}[H]
		\centering \small
		\caption{Observed yield, original crop model estimated yield, and estimated yields of five locations at each county at the end of growing season 2017 under optimal managerial decisions of different levels of risk aversion. The last column includes observed state-level yield. Percentages show the improvement of estimated yields of five locations at the end of the growing season using optimal decisions than the original estimated yield. }\label{2017ap}
		\begin{adjustbox}{width=1\textwidth}
			
			\begin{tabular}{cc|ccccc|ccccc|c|c|c}
				\clineB{1-15}{4}
				\multirow{2}{*}{\rot{State \hspace{1.5em} }}&\multirow{2}{*}{\rot{Location\hspace{1.25em}}}&\multicolumn{5}{c|}{Strategy 1}&\multicolumn{5}{c|}{Strategy 2}&\multirow{2}{*}{\rot{Strategy 3\hspace{0.9em}}}&\multirow{2}{*}{\rot{Est. Yield\hspace{0.9em}}}&\multirow{2}{*}{\rot{Obs. Yield\hspace{0.75em}}}\\
				\hhline{~~----------~~~}
				&&\rot{\hspace{0.25em}Robust}&\rot{$\alpha=0.25$}&\rot{\hspace{0.15em}$\alpha=0.5$}&\rot{$\alpha=0.75$}&\rot{Stochastic}&\rot{\hspace{0.25em}Robust}&\rot{$\alpha=0.25$}&\rot{\hspace{0.15em}$\alpha=0.5$}&\rot{$\alpha=0.75$}&\rot{Stochastic}&\multicolumn{1}{c|}{}&\multicolumn{1}{c|}{}\\
				\clineB{1-15}{3}	
				\multirow{10}{*}{Illinois}&logan-1&\makecell{230.4\\-4.9\%}&\makecell{247.7\\2.2\%}&\makecell{251\\3.6\%}&\makecell{251\\3.6\%}&\makecell{251\\3.6\%}&\makecell{251\\3.6\%}&\makecell{247.7\\2.2\%}&\makecell{247.7\\2.2\%}&\makecell{251\\3.6\%}&\makecell{251\\3.6\%}&\makecell{255.5\\5.4\%}&242.3&\multirow{10}{*}{213}\\
				\hhline{~-------------~}
				&logan-2&\makecell{236.2\\5\%}&\makecell{208.2\\-7.4\%}&\makecell{229.8\\2.2\%}&\makecell{229.8\\2.2\%}&\makecell{229.8\\2.2\%}&\makecell{208.4\\-7.4\%}&\makecell{211.9\\-5.8\%}&\makecell{218.7\\-2.8\%}&\makecell{218.7\\-2.8\%}&\makecell{229.6\\2.1\%}&\makecell{238.1\\5.9\%}&224.9\\
				\hhline{~-------------~}
				&logan-3&\makecell{215.2\\-4.5\%}&\makecell{234.2\\3.9\%}&\makecell{235\\4.3\%}&\makecell{235\\4.3\%}&\makecell{235\\4.3\%}&\makecell{234.2\\3.9\%}&\makecell{234.2\\3.9\%}&\makecell{235\\4.3\%}&\makecell{235\\4.3\%}&\makecell{235\\4.3\%}&\makecell{240.7\\6.8\%}&225.3\\
				\hhline{~-------------~}
				&logan-4&\makecell{239\\-2.5\%}&\makecell{249.7\\1.8\%}&\makecell{249.7\\1.8\%}&\makecell{249.7\\1.8\%}&\makecell{255.7\\4.3\%}&\makecell{230.6\\-5.9\%}&\makecell{249.7\\1.8\%}&\makecell{249.7\\1.8\%}&\makecell{255.8\\4.3\%}&\makecell{255.8\\4.3\%}&\makecell{264\\7.7\%}&245.2\\
				\hhline{~-------------~}
				&logan-5&\makecell{265.7\\5.3\%}&\makecell{250.8\\-0.6\%}&\makecell{260.6\\3.3\%}&\makecell{260.6\\3.3\%}&\makecell{260.6\\3.3\%}&\makecell{250.7\\-0.6\%}&\makecell{250.8\\-0.6\%}&\makecell{260.7\\3.3\%}&\makecell{260.7\\3.3\%}&\makecell{260.7\\3.3\%}&\makecell{272.6\\8.1\%}&252.2\\
				\clineB{1-15}{3}
				\multirow{10}{*}{Indiana}&greene-1&\makecell{218.2\\1.4\%}&\makecell{212.2\\-1.4\%}&\makecell{212.2\\-1.4\%}&\makecell{212.2\\-1.4\%}&\makecell{212.2\\-1.4\%}&\makecell{202.1\\-6.1\%}&\makecell{207.9\\-3.4\%}&\makecell{209.7\\-2.6\%}&\makecell{209.7\\-2.6\%}&\makecell{209.7\\-2.6\%}&\makecell{226.3\\5.1\%}&215.3&\multirow{10}{*}{161.1}\\
				\hhline{~-------------~}
				&greene-2&\makecell{208.1\\-6.3\%}&\makecell{208.1\\-6.3\%}&\makecell{211.6\\-4.8\%}&\makecell{211.6\\-4.8\%}&\makecell{208.7\\-6.1\%}&\makecell{208.1\\-6.3\%}&\makecell{208.7\\-6.1\%}&\makecell{208.7\\-6.1\%}&\makecell{208.7\\-6.1\%}&\makecell{208.7\\-6.1\%}&\makecell{225.8\\1.6\%}&222.2\\
				\hhline{~-------------~}
				&greene-3&\makecell{208.1\\-5.1\%}&\makecell{208.7\\-4.8\%}&\makecell{214.7\\-2\%}&\makecell{208.7\\-4.8\%}&\makecell{208.7\\-4.8\%}&\makecell{208.1\\-5.1\%}&\makecell{205\\-6.5\%}&\makecell{205\\-6.5\%}&\makecell{205\\-6.5\%}&\makecell{205\\-6.5\%}&\makecell{224\\2.2\%}&219.2\\
				\hhline{~-------------~}
				&greene-4&\makecell{171.9\\-26\%}&\makecell{220.4\\-5.1\%}&\makecell{223.2\\-4\%}&\makecell{223.2\\-4\%}&\makecell{223.2\\-4\%}&\makecell{207.3\\-10.8\%}&\makecell{217.3\\-6.5\%}&\makecell{217.3\\-6.5\%}&\makecell{217.3\\-6.5\%}&\makecell{217.3\\-6.5\%}&\makecell{241.8\\4.1\%}&232.3\\
				\hhline{~-------------~}
				&greene-5&\makecell{179.9\\0.2\%}&\makecell{186\\3.6\%}&\makecell{192.6\\7.3\%}&\makecell{192.6\\7.3\%}&\makecell{192.6\\7.3\%}&\makecell{186\\3.6\%}&\makecell{186\\3.6\%}&\makecell{186\\3.6\%}&\makecell{186\\3.6\%}&\makecell{197.2\\9.8\%}&\makecell{194\\8\%}&179.5\\
				\clineB{1-15}{3}
				\multirow{30}{*}{Iowa}&boone-1&\makecell{268.4\\7.8\%}&\makecell{268.3\\7.8\%}&\makecell{262.9\\5.6\%}&\makecell{249.3\\0.1\%}&\makecell{252.2\\1.3\%}&\makecell{236.6\\-5\%}&\makecell{264.9\\6.4\%}&\makecell{264.9\\6.4\%}&\makecell{264.9\\6.4\%}&\makecell{264.9\\6.4\%}&\makecell{269\\8\%}&249&\multirow{10}{*}{192.4}\\
				\hhline{~-------------~}
				&boone-2&\makecell{235\\-5.8\%}&\makecell{265.2\\6.3\%}&\makecell{265.2\\6.3\%}&\makecell{265.2\\6.3\%}&\makecell{270.7\\8.5\%}&\makecell{265.2\\6.3\%}&\makecell{265.2\\6.3\%}&\makecell{271.6\\8.8\%}&\makecell{271.6\\8.8\%}&\makecell{274.1\\9.8\%}&\makecell{271.2\\8.7\%}&249.6\\
				\hhline{~-------------~}
				&boone-3&\makecell{265.4\\9.5\%}&\makecell{265.4\\9.5\%}&\makecell{265.4\\9.5\%}&\makecell{260.4\\7.4\%}&\makecell{265.4\\9.5\%}&\makecell{265.4\\9.5\%}&\makecell{257.4\\6.2\%}&\makecell{272.9\\12.6\%}&\makecell{272.9\\12.6\%}&\makecell{272.9\\12.6\%}&\makecell{272.2\\12.3\%}&242.4\\
				\hhline{~-------------~}
				&boone-4&\makecell{243.6\\-4.6\%}&\makecell{263\\2.9\%}&\makecell{265.3\\3.9\%}&\makecell{265.3\\3.9\%}&\makecell{265.3\\3.9\%}&\makecell{249.6\\-2.3\%}&\makecell{256.9\\0.6\%}&\makecell{263.8\\3.3\%}&\makecell{265.2\\3.8\%}&\makecell{271.7\\6.3\%}&\makecell{262.2\\2.6\%}&255.5\\
				\hhline{~-------------~}
				&boone-5&\makecell{265.1\\4.3\%}&\makecell{273.3\\7.5\%}&\makecell{273.3\\7.5\%}&\makecell{273.3\\7.5\%}&\makecell{273.3\\7.5\%}&\makecell{250.4\\-1.5\%}&\makecell{267.3\\5.2\%}&\makecell{272.8\\7.4\%}&\makecell{272.8\\7.4\%}&\makecell{272.8\\7.4\%}&\makecell{272.7\\7.3\%}&254.1\\
				\clineB{2-15}{3}
				&keokuk-1&\makecell{238.1\\-5.2\%}&\makecell{251.5\\0.1\%}&\makecell{251.5\\0.1\%}&\makecell{251.5\\0.1\%}&\makecell{251.5\\0.1\%}&\makecell{239.3\\-4.7\%}&\makecell{239.3\\-4.7\%}&\makecell{239.3\\-4.7\%}&\makecell{246.1\\-2\%}&\makecell{246.1\\-2\%}&\makecell{254.2\\1.2\%}&251.2&\multirow{10}{*}{176.4}\\
				\hhline{~-------------~}
				&keokuk-2&\makecell{210.7\\-1.5\%}&\makecell{204.2\\-4.6\%}&\makecell{222.1\\3.8\%}&\makecell{222.1\\3.8\%}&\makecell{222.1\\3.8\%}&\makecell{209.2\\-2.2\%}&\makecell{209.2\\-2.2\%}&\makecell{215.4\\0.7\%}&\makecell{215.4\\0.7\%}&\makecell{215.4\\0.7\%}&\makecell{215.4\\0.7\%}&213.9\\
				\hhline{~-------------~}
				&keokuk-3&\makecell{213.6\\-2.4\%}&\makecell{213.6\\-2.4\%}&\makecell{213.6\\-2.4\%}&\makecell{218.8\\0\%}&\makecell{228.8\\4.5\%}&\makecell{213.6\\-2.4\%}&\makecell{214.9\\-1.8\%}&\makecell{214.9\\-1.8\%}&\makecell{214.9\\-1.8\%}&\makecell{218.8\\0\%}&\makecell{219.8\\0.5\%}&218.8\\
				\hhline{~-------------~}
				&keokuk-4&\makecell{211.2\\-4.5\%}&\makecell{217.8\\-1.5\%}&\makecell{218.9\\-1\%}&\makecell{218.9\\-1\%}&\makecell{218.9\\-1\%}&\makecell{211.6\\-4.3\%}&\makecell{223.8\\1.3\%}&\makecell{223.8\\1.3\%}&\makecell{228.7\\3.4\%}&\makecell{228.7\\3.4\%}&\makecell{232.2\\5\%}&221.1\\
				\hhline{~-------------~}
				&keokuk-5&\makecell{224.1\\-2.7\%}&\makecell{224.1\\-2.7\%}&\makecell{224.1\\-2.7\%}&\makecell{232.6\\1\%}&\makecell{232.6\\1\%}&\makecell{226.2\\-1.7\%}&\makecell{226.2\\-1.7\%}&\makecell{226.2\\-1.7\%}&\makecell{226.2\\-1.7\%}&\makecell{226.2\\-1.7\%}&\makecell{237\\2.9\%}&230.2\\
				\clineB{2-15}{3}
				&obrien-1&\makecell{250.5\\3.8\%}&\makecell{253.1\\4.9\%}&\makecell{253.1\\4.9\%}&\makecell{253.1\\4.9\%}&\makecell{267.5\\10.9\%}&\makecell{251\\4\%}&\makecell{251.5\\4.2\%}&\makecell{264.6\\9.7\%}&\makecell{253.1\\4.9\%}&\makecell{267.5\\10.9\%}&\makecell{271.9\\12.7\%}&241.3&\multirow{10}{*}{213.8}\\
				\hhline{~-------------~}
				&obrien-2&\makecell{255.2\\1\%}&\makecell{255.2\\1\%}&\makecell{255.2\\1\%}&\makecell{279.8\\10.8\%}&\makecell{298.8\\18.3\%}&\makecell{276\\9.3\%}&\makecell{281.9\\11.7\%}&\makecell{281.9\\11.7\%}&\makecell{281.9\\11.7\%}&\makecell{279.8\\10.8\%}&\makecell{302.5\\19.8\%}&252.5\\
				\hhline{~-------------~}
				&obrien-3&\makecell{264.6\\6.3\%}&\makecell{277.3\\11.4\%}&\makecell{286.1\\15\%}&\makecell{286.1\\15\%}&\makecell{286.1\\15\%}&\makecell{259.7\\4.4\%}&\makecell{269.9\\8.5\%}&\makecell{277.4\\11.5\%}&\makecell{277.4\\11.5\%}&\makecell{285.9\\14.9\%}&\makecell{290.9\\16.9\%}&248.8\\
				\hhline{~-------------~}
				&obrien-4&\makecell{274\\8.5\%}&\makecell{274\\8.5\%}&\makecell{274\\8.5\%}&\makecell{286.9\\13.6\%}&\makecell{286.9\\13.6\%}&\makecell{258.3\\2.2\%}&\makecell{284.1\\12.5\%}&\makecell{284.1\\12.5\%}&\makecell{284.1\\12.5\%}&\makecell{284.1\\12.5\%}&\makecell{292.6\\15.8\%}&252.6\\
				\hhline{~-------------~}
				&obrien-5&\makecell{172.7\\-14\%}&\makecell{172.7\\-14\%}&\makecell{172.7\\-14\%}&\makecell{172.7\\-14\%}&\makecell{172.7\\-14\%}&\makecell{236.3\\17.7\%}&\makecell{241.9\\20.5\%}&\makecell{241.9\\20.5\%}&\makecell{254.2\\26.6\%}&\makecell{254.2\\26.6\%}&\makecell{260.7\\29.9\%}&200.8\\
				\clineB{1-15}{4}
			\end{tabular}
		\end{adjustbox}
	\end{table}	
\end{center}

\begin{center}	
	\begin{table}[H]
		\centering \small
		\caption{Observed yield, original crop model estimated yield, and estimated yields of five locations at each county at the end of growing season 2018 under optimal managerial decisions of different levels of risk aversion. The last column includes observed state-level yield. Percentages show the improvement of estimated yields of five locations at the end of the growing season using optimal decisions than the original estimated yield. }\label{2018ap}
		\begin{adjustbox}{width=1\textwidth}
			
			\begin{tabular}{cc|ccccc|ccccc|c|c|c}
				\clineB{1-15}{4}
				\multirow{2}{*}{\rot{State \hspace{1.5em} }}&\multirow{2}{*}{\rot{Location\hspace{1.25em}}}&\multicolumn{5}{c|}{Strategy 1}&\multicolumn{5}{c|}{Strategy 2}&\multirow{2}{*}{\rot{Strategy 3\hspace{0.9em}}}&\multirow{2}{*}{\rot{Est. Yield\hspace{0.9em}}}&\multirow{2}{*}{\rot{Obs. Yield\hspace{0.75em}}}\\
				\hhline{~~----------~~~}
				&&\rot{\hspace{0.25em}Robust}&\rot{$\alpha=0.25$}&\rot{\hspace{0.15em}$\alpha=0.5$}&\rot{$\alpha=0.75$}&\rot{Stochastic}&\rot{\hspace{0.25em}Robust}&\rot{$\alpha=0.25$}&\rot{\hspace{0.15em}$\alpha=0.5$}&\rot{$\alpha=0.75$}&\rot{Stochastic}&\multicolumn{1}{c|}{}&\multicolumn{1}{c|}{}\\
				\clineB{1-15}{3}
				\multirow{10}{*}{Illinois}&logan-1&\makecell{204.3\\-3.9\%}&\makecell{204.3\\-3.9\%}&\makecell{220.7\\3.9\%}&\makecell{220.7\\3.9\%}&\makecell{220.7\\3.9\%}&\makecell{209.8\\-1.3\%}&\makecell{223.4\\5.1\%}&\makecell{223.4\\5.1\%}&\makecell{223.4\\5.1\%}&\makecell{223.4\\5.1\%}&\makecell{234.5\\10.3\%}&212.5&\multirow{10}{*}{236.2}\\
				\hhline{~-------------~}
				&logan-2&\makecell{201.8\\-6.3\%}&\makecell{201.8\\-6.3\%}&\makecell{211.2\\-1.9\%}&\makecell{211.2\\-1.9\%}&\makecell{210.8\\-2.1\%}&\makecell{213.3\\-0.9\%}&\makecell{215.6\\0.2\%}&\makecell{215.6\\0.2\%}&\makecell{218.2\\1.4\%}&\makecell{217.3\\1\%}&\makecell{221.6\\3\%}&215.2\\
				\hhline{~-------------~}
				&logan-3&\makecell{208.4\\-1.3\%}&\makecell{216.4\\2.5\%}&\makecell{216.4\\2.5\%}&\makecell{216.4\\2.5\%}&\makecell{216.4\\2.5\%}&\makecell{212.7\\0.7\%}&\makecell{217.9\\3.2\%}&\makecell{217.9\\3.2\%}&\makecell{220.5\\4.4\%}&\makecell{220.5\\4.4\%}&\makecell{226.1\\7.1\%}&211.2\\
				\hhline{~-------------~}
				&logan-4&
				\makecell{235.1\\5.3\%}&\makecell{226.3\\1.4\%}&\makecell{229.4\\2.7\%}&\makecell{229.4\\2.7\%}&\makecell{229.4\\2.7\%}&\makecell{220.7\\-1.2\%}&\makecell{234.7\\5.1\%}&\makecell{238\\6.6\%}&\makecell{238\\6.6\%}&\makecell{238\\6.6\%}&\makecell{242.8\\8.7\%}&223.3\\
				\hhline{~-------------~}
				&logan-5&\makecell{223.8\\0.1\%}&\makecell{226.9\\1.6\%}&\makecell{229.6\\2.8\%}&\makecell{229.6\\2.8\%}&\makecell{229.6\\2.8\%}&\makecell{232.4\\4\%}&\makecell{235.4\\5.3\%}&\makecell{238.7\\6.8\%}&\makecell{238.7\\6.8\%}&\makecell{238.7\\6.8\%}&\makecell{243.8\\9.1\%}&223.5\\
				\clineB{1-15}{3}
				\multirow{10}{*}{Indiana}&greene-1&\makecell{171.5\\-15.2\%}&\makecell{214.5\\6.1\%}&\makecell{214.5\\6.1\%}&\makecell{214.5\\6.1\%}&\makecell{214.5\\6.1\%}&\makecell{143.5\\-29.1\%}&\makecell{210.2\\4\%}&\makecell{214.4\\6\%}&\makecell{214.4\\6\%}&\makecell{214.4\\6\%}&\makecell{217.6\\7.6\%}&202.2&\multirow{10}{*}{178.9}\\
				\hhline{~-------------~}
				&greene-2&\makecell{173.3\\-10.1\%}&\makecell{173.3\\-10.1\%}&\makecell{173.3\\-10.1\%}&\makecell{209.9\\8.9\%}&\makecell{209.9\\8.9\%}&\makecell{173.4\\-10.1\%}&\makecell{173.3\\-10.1\%}&\makecell{173.4\\-10.1\%}&\makecell{209.9\\8.9\%}&\makecell{209.9\\8.9\%}&\makecell{214.9\\11.4\%}&192.8\\
				\hhline{~-------------~}
				&greene-3&\makecell{143.9\\-21.7\%}&\makecell{205.1\\11.6\%}&\makecell{205.1\\11.6\%}&\makecell{205.1\\11.6\%}&\makecell{205.1\\11.6\%}&\makecell{140.3\\-23.7\%}&\makecell{209.1\\13.8\%}&\makecell{209.1\\13.8\%}&\makecell{205.1\\11.6\%}&\makecell{205.1\\11.6\%}&\makecell{209.1\\13.8\%}&183.8\\
				\hhline{~-------------~}
				&greene-4&\makecell{126.9\\-38.6\%}&\makecell{212.8\\2.9\%}&\makecell{217.8\\5.3\%}&\makecell{217.8\\5.3\%}&\makecell{217.8\\5.3\%}&\makecell{126.9\\-38.6\%}&\makecell{212.8\\2.9\%}&\makecell{217.8\\5.3\%}&\makecell{217.8\\5.3\%}&\makecell{217.8\\5.3\%}&\makecell{222\\7.4\%}&206.8\\
				\hhline{~-------------~}
				&greene-5&\makecell{122\\-30.5\%}&\makecell{192\\9.4\%}&\makecell{179.6\\2.4\%}&\makecell{168.3\\-4.1\%}&\makecell{168.3\\-4.1\%}&\makecell{126.6\\-27.9\%}&\makecell{192\\9.4\%}&\makecell{196.2\\11.8\%}&\makecell{196.2\\11.8\%}&\makecell{196.2\\11.8\%}&\makecell{198.7\\13.2\%}&175.5\\
				\clineB{1-15}{3}
				\multirow{30}{*}{Iowa}&boone-1&\makecell{236.4\\9.1\%}&\makecell{230.5\\6.4\%}&\makecell{230.5\\6.4\%}&\makecell{230.5\\6.4\%}&\makecell{230.5\\6.4\%}&\makecell{236.4\\9.1\%}&\makecell{236.4\\9.1\%}&\makecell{233.6\\7.8\%}&\makecell{233.6\\7.8\%}&\makecell{230.5\\6.4\%}&\makecell{242.1\\11.7\%}&216.7&\multirow{10}{*}{193.7}\\
				\hhline{~-------------~}
				&boone-2&\makecell{220.9\\1.9\%}&\makecell{220\\1.5\%}&\makecell{223.4\\3\%}&\makecell{223.4\\3\%}&\makecell{223.4\\3\%}&\makecell{220.9\\1.9\%}&\makecell{223.4\\3\%}&\makecell{223.4\\3\%}&\makecell{223.4\\3\%}&\makecell{223.4\\3\%}&\makecell{238.1\\9.8\%}&216.8\\
				\hhline{~-------------~}
				&boone-3&\makecell{211.6\\-0.6\%}&\makecell{220.8\\3.8\%}&\makecell{220.8\\3.8\%}&\makecell{223.6\\5.1\%}&\makecell{223.6\\5.1\%}&\makecell{217\\2\%}&\makecell{217.1\\2\%}&\makecell{222.4\\4.5\%}&\makecell{222.4\\4.5\%}&\makecell{222.4\\4.5\%}&\makecell{235.2\\10.5\%}&212.8\\
				\hhline{~-------------~}
				&boone-4&\makecell{216.1\\-0.2\%}&\makecell{225.5\\4.1\%}&\makecell{226.6\\4.7\%}&\makecell{226.6\\4.7\%}&\makecell{226.6\\4.7\%}&\makecell{220.2\\1.7\%}&\makecell{225.5\\4.1\%}&\makecell{225.5\\4.1\%}&\makecell{226.6\\4.7\%}&\makecell{226.6\\4.7\%}&\makecell{235\\8.5\%}&216.6\\
				\hhline{~-------------~}
				&boone-5&\makecell{225.3\\4.1\%}&\makecell{221.9\\2.5\%}&\makecell{228.6\\5.6\%}&\makecell{228.6\\5.6\%}&\makecell{228.6\\5.6\%}&\makecell{221.9\\2.5\%}&\makecell{221.9\\2.5\%}&\makecell{224.6\\3.7\%}&\makecell{224.6\\3.7\%}&\makecell{224.5\\3.7\%}&\makecell{246.9\\14\%}&216.5\\
				\clineB{2-15}{3}
				&keokuk-1&\makecell{224\\3.4\%}&\makecell{225.6\\4.1\%}&\makecell{224\\3.4\%}&\makecell{224\\3.4\%}&\makecell{224.9\\3.8\%}&\makecell{224.2\\3.4\%}&\makecell{222\\2.5\%}&\makecell{225.6\\4.1\%}&\makecell{225.6\\4.1\%}&\makecell{225.6\\4.1\%}&\makecell{238\\9.8\%}&216.7&\multirow{10}{*}{203.8}\\
				\hhline{~-------------~}
				&keokuk-2&\makecell{207.4\\3.7\%}&\makecell{203.6\\1.9\%}&\makecell{205\\2.6\%}&\makecell{205\\2.6\%}&\makecell{205\\2.6\%}&\makecell{214.3\\7.2\%}&\makecell{214.3\\7.2\%}&\makecell{214.3\\7.2\%}&\makecell{214.3\\7.2\%}&\makecell{213.3\\6.7\%}&\makecell{219.4\\9.7\%}&199.9\\
				\hhline{~-------------~}
				&keokuk-3&\makecell{220.3\\6.4\%}&\makecell{211.8\\2.3\%}&\makecell{217.4\\5\%}&\makecell{217.4\\5\%}&\makecell{209.7\\1.3\%}&\makecell{218.1\\5.4\%}&\makecell{214.9\\3.8\%}&\makecell{214.9\\3.8\%}&\makecell{214.9\\3.8\%}&\makecell{214.9\\3.8\%}&\makecell{226.5\\9.4\%}&207\\
				\hhline{~-------------~}
				&keokuk-4&\makecell{209.4\\2.4\%}&\makecell{201.8\\-1.3\%}&\makecell{203.1\\-0.7\%}&\makecell{203.1\\-0.7\%}&\makecell{203.1\\-0.7\%}&\makecell{205.3\\0.4\%}&\makecell{213.5\\4.4\%}&\makecell{212.9\\4.1\%}&\makecell{209.2\\2.3\%}&\makecell{209.2\\2.3\%}&\makecell{224.8\\10\%}&204.5\\
				\hhline{~-------------~}
				&keokuk-5&\makecell{216\\4\%}&\makecell{216\\4\%}&\makecell{216\\4\%}&\makecell{216\\4\%}&\makecell{211.2\\1.7\%}&\makecell{204.9\\-1.4\%}&\makecell{216.2\\4.1\%}&\makecell{216.2\\4.1\%}&\makecell{216.2\\4.1\%}&\makecell{210.4\\1.3\%}&\makecell{225.1\\8.4\%}&207.7\\
				\clineB{2-15}{3}
				&obrien-1&\makecell{240.8\\3.4\%}&\makecell{240.2\\3.1\%}&\makecell{240.2\\3.1\%}&\makecell{240.2\\3.1\%}&\makecell{235.6\\1.2\%}&\makecell{242.9\\4.3\%}&\makecell{240.4\\3.2\%}&\makecell{240.4\\3.2\%}&\makecell{240.4\\3.2\%}&\makecell{240.4\\3.2\%}&\makecell{252.1\\8.2\%}&232.9&\multirow{10}{*}{201.8}\\
				\hhline{~-------------~}
				&obrien-2&\makecell{243.8\\4.3\%}&\makecell{242.5\\3.8\%}&\makecell{242.5\\3.8\%}&\makecell{241.7\\3.4\%}&\makecell{241.7\\3.4\%}&\makecell{243.8\\4.3\%}&\makecell{242.5\\3.8\%}&\makecell{242.5\\3.8\%}&\makecell{242.4\\3.7\%}&\makecell{241.7\\3.4\%}&\makecell{251.7\\7.8\%}&233.6\\
				\hhline{~-------------~}
				&obrien-3&\makecell{242.4\\2.3\%}&\makecell{239.6\\1.2\%}&\makecell{239.6\\1.2\%}&\makecell{239.6\\1.2\%}&\makecell{234.3\\-1.1\%}&\makecell{240.1\\1.4\%}&\makecell{239.6\\1.2\%}&\makecell{239.6\\1.1\%}&\makecell{239.6\\1.1\%}&\makecell{239.6\\1.1\%}&\makecell{249\\5.1\%}&236.9\\
				\hhline{~-------------~}
				&obrien-4&\makecell{241.5\\3.6\%}&\makecell{239.4\\2.7\%}&\makecell{239.4\\2.7\%}&\makecell{234.2\\0.5\%}&\makecell{234.2\\0.5\%}&\makecell{245.8\\5.5\%}&\makecell{239.4\\2.7\%}&\makecell{239.4\\2.7\%}&\makecell{239.4\\2.7\%}&\makecell{234.2\\0.5\%}&\makecell{249.1\\6.9\%}&233.1\\
				\hhline{~-------------~}
				&obrien-5&\makecell{231\\4\%}&\makecell{239.8\\8\%}&\makecell{237.2\\6.8\%}&\makecell{237.2\\6.8\%}&\makecell{237\\6.7\%}&\makecell{231\\4\%}&\makecell{239.8\\8\%}&\makecell{247.8\\11.6\%}&\makecell{247.8\\11.6\%}&\makecell{247.8\\11.6\%}&\makecell{253\\13.9\%}&222.1\\
				\clineB{1-15}{4}
			\end{tabular}
		\end{adjustbox}
	\end{table}	
\end{center}

	\subsection*{Parallel Bayesian Optimization}

This section proposes an extended version of BO \cite{mockus2012bayesian} to improve the time-dependent parameter calibration structure in the proposed framework (see Figure \ref{window}). We briefly review the classic BO method and the key statistical and optimization methods on which it relies in the Appendix 1. We use BO because it is one of the most successful approaches in optimizing calibration problems \cite{jones2001taxonomy}. BO suffers from some limitations, and thus we modify this algorithm. First, BO has a high sensitivity to the type of acquisition function as an objective function and the nature of the gaussian procedure (GP), including kernel types and kernel hyperparameters \cite{ snoek2012practical}. Second, BO loses its efficiency in exploration and exploitation as the input dimensions increase. Third, the acquisition function is often difficult to optimize, and its performance depends on an optimizer to search the surface.

We develop a new parallel optimization scheme for the optimization part of the time-dependent parameter calibration framework, which addresses the BO method's drawbacks. In this approach, instead of allocating all searching budget (the number of iteration) to one BO model, we divide the searching budget into several BO algorithms to search a solution space in parallel. To address BO's sensitivity to the type of acquisition functions, each parallel BO has a specific type of acquisition function. Besides, each BO applies a particular kernel type and hyper-parameters for constructing the GP distribution. When BO models determine the next sample point by maximizing its equations function, they share their findings to help each other to acquire more knowledge about the solution space. The structural differences between BO and PBO are illustrated in Figure \ref{BO_PBO}. The PBO algorithm is defined in Algorithm \ref{PBO}.

%	\begin{algorithm}[H]
%			\caption{\texttt{Parallel Bayesian optimization}} \label{PBO}
%			\begin{algorithmic}[1]
%					\STATE \textbf{Input:} Data set $D=\{ (x_0, y_0)\}$, $N$ as number BO models are run in parallel and $T$ as maximum iteration.
%					\STATE \textbf{Output:} A local optimal solution $x^* \in \mathbb{R}^{1 \times p}$.
%					\STATE \textbf{Step 0:} Set $K_n$, $\sigma$, and $u_n$ as type of kernel, set of kernel parameters, and acquisition function of $n$th BO for all $n \in \{1,...,N\}$. Set incumbent solution $x^*=x_0$.
%					\FOR{$t=1$ to $T$}
%					\FOR{$n=1$ to $N$}
%					\STATE \textbf{Step 1:} Use GP to update the posterior probability $\hat{f}_n$ and construct acquisition function $u_n$.
%					\STATE \textbf{Step 2:} Find $x_t^n$ by optimizing the acquisition function $u_n$ over function $f$: $x_t^n = \argmax_{x}u_n(x|D)$.
%					\STATE \textbf{Step 3:} Sample $x_t^n$ to calculate the objective function $y_t^n=f(x_t^n)$ and augment the data $D = \{D,(x_t^n,y_t^n)\}$.
%					\STATE \textbf{Step 4:} Update incumbent solution $x^*  = x_{\argmax (y)}$.
%					\ENDFOR
%					\ENDFOR
%				\end{algorithmic}
%		\end{algorithm}

	\begin{algorithm}[H]
			\caption{\texttt{Parallel Bayesian optimization}} \label{PBO}
			\begin{algorithmic}[1]
					\STATE \textbf{Input:} Data set $D=\{ (x_0, y_0)\}$, $N$ as number BO models are run in parallel and $T$ as maximum iteration, $I_s$ as scenario $s \in S$, where $S$ denotes set of scenarios.
					\STATE \textbf{Output:} A local optimal solution $x^* \in \mathbb{R}^{1 \times p}$.
					\STATE \textbf{Step 0:} Set $K_n$, $\sigma$, and $u_n$ as type of kernel, set of kernel parameters, and acquisition function of $n$th BO for all $n \in \{1,...,N\}$. Set incumbent solution $x^*=x_0$.
					\FOR{$t=1$ to $T$}
					\FOR{$n=1$ to $N$}
					\STATE \textbf{Step 1:} Use GP to update the posterior probability $\hat{f}_n$ and construct acquisition function $u_n$.
					\STATE \textbf{Step 2:} Find $x_t^n$ by optimizing the acquisition function $u_n$ over function $f$: $x_t^n = \argmax_{x}u_n(x|D)$.
					\STATE \textbf{Step 3:} For each scenario $I_s$, call crop model by modifying its file according to each scenario to evalaue the performance of solution $x_t^n$ for scenario $s$ as $f(x_t^n,I_s)$. 
					\STATE \textbf{Step 4:} Calculate the objective function of solution $x_t^n$ for stochastic optimization as $y_t^n=\sum_{s \in S}{f(x_t^n,I_s)}/|S|$ and for robust optimization as $y_t^n=min_{s \in S}\{f(x_t^n,I_s)\}$ or $y_t^n=max_{s \in S}\{f(x_t^n,I_s)\}$.
					\STATE \textbf{Step 5:} Augment the data $D = \{D,(x_t^n,y_t^n)\}$ and update incumbent solution $x^*  = x_{\argmax (y)}$.
					\ENDFOR
					\ENDFOR
				\end{algorithmic}
		\end{algorithm}

	We provide an additional explanation about the algorithm as follows.
	
	\begin{itemize}
			
			\item \textbf{Hyperparameters:} Parameter $N$ should reflect the number of parallel BO models such that any increase in $N$ will greatly expand the searching power to search through more enumeration. The maximum number of iterations $T$ indicates the stop criteria of the algorithm. The maximum number of iterations can be replaced with a threshold so that we terminate the algorithm if the difference between the current value of the best new points and the incumbent solution is less than a predefined threshold.

			\item \textbf{Step 0:}  In this step, we set one type of acquisition function, kernel, and kernel parameters for constructing the GP of each BO model.
			
			\item \textbf{Step 1:} GP updates the posterior probability, and constructs acquisition function in this step.
			
			\item \textbf{Step 2:} We use the limited-memory quasi-Newton algorithm for bound-constrained optimization $\text{L-BFGS-B}$ \cite{zhu1997algorithm,byrd1995limited} in Python to optimize the acquisition function with a set of random starting points. We run $\text{L-BFGS-B}$ several times with different arbitrary starting points to improve the efficiency and overcome the difficulty of optimizing the acquisition function. We select the next sample point as the best optimal solution between all found solutions of $\text{L-BFGS-B}$ algorithm with different starting points.
			
			\item \textbf{Step 3:} The APSIM is run to evaluate a new parameter combination from each BO model in this step. Then, the BO models share their findings to expand their knowledge about posterior distribution such that the BO model applies GP in the next iteration to update its posterior distribution and construct the acquisition function.
			
			\item \textbf{Step 4:} The incumbent solution is updated as the best combinations of parameters evaluated by APSIM in this step.
		\end{itemize}

\end{document}